\input ./style/arxiv-general.cfg
\documentclass[aos,MSNbibl,nameyear,dvips]{arximspdf}
\makeatletter
   \@ifpackageloaded{graphicx}{}{\usepackage{graphicx}}
\makeatother
\usepackage{mathbh}
\usepackage[ruled,linesnumbered]{algorithm2e}

\doi{10.1214/15-AOS1321}% Updated by VTEXPTS2LaTeX.exe, 04.03.2015
\volume{43}
\issue{4}
\pubyear{2015}
\firstpage{1716}
\lastpage{1741}
\docsubty{FLA}

\makeatletter
\let\my@algocf@latexcaption\algocf@latexcaption
\let\my@addcontentsline\addcontentsline
\long\def\algocf@latexcaption#1[#2]#3{%
\def\addcontentsline##1##2##3{}%
\my@algocf@latexcaption{#1}[#2]{#3}%
\global\let\addcontentsline\my@addcontentsline%
}
\newcommand{\rrvert}{\vert}
\newcommand{\llvert}{\vert}
\newproclaim{rem}{Remark}
\newtheorem{theorem}{Theorem}
\newtheorem{lemme}{Lemma}

\newtheorem{proposition}{Proposition}

\makeatother

\begin{document}
\begin{frontmatter}

\title{Consistency of random forests\thanksref{T1}}
\runtitle{Consistency of random forests}

\begin{aug}
% Corresponding author: Erwan Scornet - erwan.scornet@upmc.fr% Updated by VTEXPTS2LaTeX.exe, 09.03.2015 08:21
%by VTEXPTS2LaTeX.exe, 04.03.2015 08:24
\author[A]{\fnms{Erwan}~\snm{Scornet}\corref{}\thanksref{m1}\ead[label=e1]{erwan.scornet@upmc.fr}},
\author[A]{\fnms{G\'erard}~\snm{Biau}\thanksref{m1}\ead[label=e2]{gerard.biau@upmc.fr}}
\and
\author[B]{\fnms{Jean-Philippe}~\snm{Vert}\thanksref{m2}\ead[label=e3]{jean-philippe.vert@mines-paristech.fr}}
\runauthor{E. Scornet, G. Biau and J.-P. Vert}
\affiliation{Sorbonne Universit\'es\thanksmark{m1} and MINES
ParisTech, PSL-Research University\thanksmark{m2}}
%\dedicated{}
\address[A]{E. Scornet\\
G. Biau\\
Sorbonne Universit\'es\\
UPMC Univ Paris 06\\
Paris F-75005\\
France
\\
\printead{e1}\\
\phantom{E-mail:\ }\printead*{e2}}
\address[B]{J.-P. Vert\\
MINES ParisTech, PSL-Research University\\
CBIO-Centre for Computational Biology\\
Fontainebleau F-77300\\
France\\
and\\
Institut Curie\\
Paris F-75248\\
France\\
and\\
U900, INSERM\\
Paris F-75248\\
France\\
\printead{e3}}
\end{aug}
\thankstext{T1}{Supported by the European Research Council
[SMAC-ERC-280032].}

% HISTORY:
%
\received{\smonth{5} \syear{2014}}% Updated by VTEXPTS2LaTeX.exe,
%04.03.2015 08:24
%
\revised{\smonth{2} \syear{2015}}% Updated by VTEXPTS2LaTeX.exe,
%04.03.2015 08:24

% ABSTRACT
%
\begin{abstract}
Random forests are a learning algorithm proposed by
Breiman [\textit{Mach. Learn.} \textbf{45} (2001) 5--32] that combines
several randomized decision trees and aggregates their predictions by
averaging. Despite its wide usage and outstanding practical
performance, little is known about the mathematical properties of the
procedure. This disparity between theory and practice originates in the
difficulty to simultaneously analyze both the randomization process and
the highly data-dependent tree structure.
In the present paper, we take a step forward in forest exploration by
proving a consistency result for Breiman's [\textit{Mach. Learn.}  \textbf{45}
(2001) 5--32] original
algorithm in the context of additive regression models.
Our analysis also sheds an interesting light on how random forests can
nicely adapt to sparsity.
\end{abstract}

% KEYWORDS
\end{frontmatter}

\section{Introduction}\label{sec1}

Random forests are an ensemble learning method for classification and
regression that constructs a number of randomized decision trees during
the training phase and predicts by averaging the results. Since its
publication in the seminal paper of \citet{RF}, the procedure has
become a major data analysis tool, that performs well in practice in
comparison with many standard methods. What has greatly contributed to
the popularity of forests is the fact that they can be applied to a
wide range of prediction problems and have few parameters to tune.
Aside from being simple to use, the method is generally recognized for
its accuracy and its ability to deal with small sample sizes,
high-dimensional feature spaces and complex data structures.
The random forest methodology has been successfully involved in many
practical problems, including air quality prediction (winning code of
the EMC data science global hackathon in 2012, see \surl{http://www.kaggle.com/c/dsg-hackathon}),
chemoinformatics [\citet{SvLiToCuShFe03}], ecology [\citet
{PrIvLi06,CuEdBe07}], 3D object
recognition [\citet{ShFiCoShFiMoKiBl11}] and bioinformatics
[\citet{microarray}],
just to name a few. In addition, many variations on
the original algorithm have been proposed to improve the calculation
time while maintaining good prediction accuracy; see, for example,
\citet{extra,ERF}. Breiman's forests have also been extended to
quantile estimation [\citet{Me06}], survival analysis [\citet{RSF}] and
ranking prediction [\citet{RankingF}].

On the theoretical side, the story is less conclusive, and regardless
of their extensive use in practical settings, little is known about the
mathematical properties of random forests. To date, most studies have
concentrated on isolated parts or simplified versions of the procedure.
The most celebrated theoretical result is that of \citet{RF}, which
offers an upper bound on the generalization error of forests in terms
of correlation and strength of the individual trees. This was followed
by a technical note [\citet{Br04}] that focuses on a stylized
version of the original algorithm. A critical step was subsequently
taken by \citet{LiJe06}, who established lower bounds for nonadaptive
forests (i.e., independent of the training set). They also highlighted
an interesting connection between random forests and a particular class
of nearest neighbor predictors that was further worked out by
\citet{BiDe10}. In recent years, various theoretical studies
[e.g., \citet{consistency,IsKo10,analysis,Ge12,RLT}] have been performed,
analyzing consistency of simplified models, and moving ever closer to
practice. Recent attempts toward narrowing the gap between theory and
practice are by \citet{DeMaFr13}, who proves the first consistency
result for online random forests, and by \citet{Wa14} and
\citet{MeHo14a} who study the asymptotic sampling distribution of forests.

The difficulty in properly analyzing random forests can be explained by
the black-box nature of the procedure, which is actually a subtle
combination of different components. Among the forest essential
ingredients, both bagging [\citet{Br96}] and the classification and
regression trees (CART)-split criterion [\citet{CART}] play a critical
role. Bagging (a contraction of bootstrap-aggregating) is a general
aggregation scheme which proceeds by generating subsamples from the
original data set, constructing a predictor from each resample and
deciding by averaging. It is one of the most effective computationally
intensive procedures to improve on unstable estimates, especially for
large, high-dimensional data sets where finding a good model in one
step is impossible because of the complexity and scale of the problem
[\citet{Analyzingbagging,KlTaSaJo12,WaHaEf13}]. The CART-split
selection originated from the most influential CART algorithm of
\citet{CART}, and is used in the construction of the individual
trees to
choose the best cuts perpendicular to the axes. At each node of each
tree, the best cut is selected by optimizing the CART-split criterion,
based on the notion of Gini impurity (classification) and prediction
squared error (regression).

Yet, while bagging and the CART-splitting scheme play a key role in the
random forest mechanism, both are difficult to analyze, thereby
explaining why theoretical studies have, thus far, considered
simplified versions of the original procedure. This is often done by
simply ignoring the bagging step and by replacing the CART-split
selection with a more elementary cut protocol.
Besides, in Breiman's forests, each leaf (i.e., a terminal node) of the
individual trees contains a fixed pre-specified number of observations
(this parameter, called \texttt{nodesize} in the R package \texttt
{randomForests}, is usually chosen between $1$ and $5$). There is also
an extra parameter in the algorithm which allows one to control the
total number of leaves (this parameter is called \texttt{maxnode} in
the R package and has, by default, no effect on the procedure). The
combination of these various components makes the algorithm difficult
to analyze with rigorous mathematics. As a matter of fact, most authors
focus on simplified, data-independent procedures, thus creating a gap
between theory and practice.

Motivated by the above discussion, we study in the present paper some
asymptotic properties of Breiman's (\citeyear{RF}) algorithm in the
context of additive regression models. We prove the $\mathbb{L}^2$
consistency of random forests, which gives a first basic theoretical
guarantee of efficiency for this algorithm.
To our knowledge, this is the first consistency result for Breiman's
(\citeyear{RF}) original procedure.
% since most of the previous studies focus on data-independent
%splitting criteria and cells containing a number of points growing to
%infinity with the sample size.
Our approach rests upon a detailed analysis of the behavior of the
cells generated by CART-split selection as the sample size grows.
It turns out that a good control of the regression function variation
inside each cell, together with a proper choice of the total number of
leaves (Theorem~\ref{theoremconsistencysemidevelopedBRF}) or a
proper choice of the subsampling rate (Theorem~\ref
{Theoremconsistanceforetbreiman}) are sufficient to ensure the
forest consistency in a $\mathbb{L}^2$ sense.
Also, our analysis shows that random forests can adapt to a sparse
framework, when the ambient dimension $p$ is large (independent of
$n$), but only a smaller number of coordinates carry out information.

The paper is organized as follows.
In Section~\ref{sec2}, we introduce some notation and describe the random
forest method. The main asymptotic results are presented in
Section~\ref{sec3}
and further discussed in Section~\ref{sec4}. Section~\ref{sec5} is
devoted to the
main proofs, and technical results are gathered in the supplemental
article [\citet{ScBiVE15Supp}].

\section{Random forests}\label{sec2}

The general framework is $\mathbb{L}^2$ regression estimation, in
which an input random vector $\mathbf{X}\in[0,1]^p$ is observed, and the
goal is to predict the square integrable random response $Y\in\mathbb
R$ by estimating the regression function $m(\mathbf{x})=\mathbb
E[Y|\mathbf{X}=\mathbf{x}]$. To this end, we assume given a training sample
$\mathcal D_n=(\mathbf{X}_1,Y_1),\ldots, (\mathbf{X}_n,Y_n)$ of
$[0,1]^p\times
\mathbb{R}$-valued independent random variables distributed as the
independent prototype pair $(\mathbf{X}, Y)$. The objective is to use the
data set $\mathcal D_n$ to construct an estimate $m_n\dvtx  [0,1]^p \to
\mathbb R$ of the function $m$. In this respect, we say that a
regression function estimate $m_n$ is $\mathbb{L}^2$ consistent if
$\mathbb{E} [m_n(\mathbf{X})-m(\mathbf{X})]^2 \to0$ as $n \to
\infty$ (where the
expectation is over $\mathbf{X}$ and $\mathcal{D}_n$).

A random forest is a predictor consisting of a collection of $M$
randomized regression trees. For the $j$th tree in the family, the
predicted value at the query point $\mathbf{x}$ is denoted by
$m_n(\mathbf{x}; \Theta_j,\mathcal{D}_n)$, where $\Theta_1,\ldots,\Theta_M$ are independent random variables, distributed as a generic
random variable $\Theta$ and independent of $\mathcal{D}_n$.
In practice, this variable is used to resample the training set prior
to the growing of individual trees and to select the successive
candidate directions for splitting. The trees are combined to form the
(finite) forest estimate
\begin{equation}
\label{finiteforest}
m_{M,n}(\mathbf{x}; \Theta_1,\ldots,
\Theta_M, \mathcal {D}_n)=\frac{1}{M}\sum
_{j=1}^M m_n(\mathbf{x};
\Theta_j,\mathcal D_n).
\end{equation}
Since in practice we can choose $M$ as large as possible, we study in
this paper the property of the infinite forest estimate obtained as the
limit of (\ref{finiteforest}) when the number of trees $M$ grows to
infinity as follows:
\begin{eqnarray*}
&& m_{n}(\mathbf{x}; \mathcal{D}_n)=\mathbb
E_{\Theta} \bigl[m_n(\mathbf{x};\Theta,\mathcal
D_n) \bigr],
\end{eqnarray*}
where $\mathbb E_{\Theta}$ denotes expectation with respect to the
random parameter $\Theta$, conditional on $\mathcal D_n$. This
operation is justified by the law of large numbers, which asserts that,
almost surely, conditional on $\mathcal{D}_n$,
\[
\lim_{M \to\infty} m_{n,M}(\mathbf{x}; \Theta_1,\ldots, \Theta_M, \mathcal{D}_n) = m_n(
\mathbf{x}; \mathcal{D}_n);
\]
see, for example, \citet{Sc14,RF} for details.
In the sequel, to lighten notation, we will simply write $m_n(\mathbf
{x})$ instead of $m_n(\mathbf{x}$; $ \mathcal{D}_n)$.

{\DecMargin{-0.3em}
\begin{algorithm}
\caption{Breiman's random forest predicted value at $\mathbf{x}$}\label{al1}
 \KwIn{Training set $\mathcal{D}_n$, number of trees $M>0$, $m_{\mathrm{try}} \in \{1, \ldots, p \}$, $a_n  \in \{1, \ldots, n \}$,  $t_n \in \{1, \ldots, a_n \}$, and $\mathbf{x} \in [0,1]^p$.}

 \KwOut{Prediction of the random forest at $\mathbf{x}$.}

 \For{$j=1, \ldots, M$}{
 Select $a_n$ points, without replacement, uniformly in $\mathcal{D}_n$.

 Set $\mathcal{P}_0 = \{[0,1]^p\}$ the partition associated with the root of the tree.

For all $1 \leq \ell \leq a_n$, set $\mathcal{P}_{\ell} = \varnothing$.

Set $n_{\mathrm{nodes}} = 1$ and $\mathrm{level} = 0$.

\While{ $n_{\mathrm{nodes}}<t_n$}{
    \eIf{$\mathcal{P}_{\mathrm{level}} = \varnothing$}{$\mathrm{level} = \mathrm{level} +1$}
    {Let $A$ be the first element in $\mathcal{P}_{\mathrm{level}}$.

        \eIf{$A$ contains exactly one point}{
                    $\mathcal{P}_{\mathrm{level}} \leftarrow \mathcal{P}_{\mathrm{level}} \setminus
                    \{A\}$\\
                    $\mathcal{P}_{\mathrm{level}+1} \leftarrow \mathcal{P}_{\mathrm{level}+1} \cup \{A\}$
        }
                {
        Select uniformly, without replacement, a subset $\mathcal{M}_{\mathrm{try}} \subset \{1,\ldots,p\}$ of cardinality
        $m_{\mathrm{try}}$.\\
        Select the best split in $A$ by optimizing the CART-split criterion along the coordinates in                    $\mathcal{M}_{\mathrm{try}}$ (\textit{see details
        below}).\\
        Cut the cell $A$ according to the best split. Call $A_L$ and $A_R$ the two resulting
        cell.\\
        $\mathcal{P}_{\mathrm{level}} \leftarrow
        \mathcal{P}_{\mathrm{level}}\setminus\{A\}$\\
        $\mathcal{P}_{\mathrm{level}+1} \leftarrow \mathcal{P}_{\mathrm{level}+1}\cup \{A_L\} \cup
        \{A_R\}$\\
        $n_{\mathrm{nodes}} = n_{\mathrm{nodes}} +1$
                }}
}
 Compute the predicted value $m_n(\mathbf{x}; \Theta_j, \mathcal{D}_n)$ at $\mathbf{x}$ equal to the  average of the $Y_i$'s falling in the
 cell of $\mathbf{x}$ in partition $\mathcal{P}_{\mathrm{level}} \cup \mathcal{P}_{\mathrm{level}+1}$.
 }
 Compute the random forest estimate $m_{M,n}(\mathbf{x}; \Theta_1, \ldots, \Theta_M, \mathcal{D}_n)$ at the query point $\mathbf{x}$ according to  (\ref{finiteforest}).
\end{algorithm}}

In Breiman's (\citeyear{RF}) original forests, each node of a single
tree is associated with a hyper-rectangular cell. At each step of the
tree construction, the collection of cells forms a partition of
$[0,1]^p$. The root of the tree is $[0,1]^p$ itself, and each tree is
grown as
explained in  Algorithm~\ref{al1}.

This algorithm has three parameters:
\begin{enumerate}[(3)]
\item[(1)] $m_{\mathrm{try}}\in\{1,\ldots, p \}$, which is
the number of pre-selected directions for splitting;
\item[(2)] $a_n \in\{1,\ldots, n \}$, which is the number of
sampled data points in each tree;
\item[(3)] $t_n \in\{1,\ldots, a_n \}$, which is the number of
leaves in each tree.
\end{enumerate}
By default, in the original procedure, the parameter $m_{\mbox{\tiny
try}}$ is set to $p/3$, $a_n$ is set to~$n$ (resampling is done with
replacement) and $t_n =a_n$.
However, in our approach, resampling is done without replacement and
the parameters $a_n$, and $t_n$ can be different from their default values.

In words, the algorithm works by growing $M$ different trees as
follows. For each tree, $a_n$ data points are drawn at random without
replacement from the original data set; then, at each cell of every
tree, a split is chosen by maximizing the CART-criterion (see below);
finally, the construction of every tree is stopped when the total
number of cells in the tree reaches the value $t_n$ (therefore, each
cell contains exactly one point in the case $t_n=a_n$).

We note that the resampling step in  Algorithm~\ref{al1} (line $2$) is
done by choosing $a_n$ out of $n$ points (with $a_n \leq n$) without
replacement. This is slightly different from the original algorithm,
where resampling is done by bootstrapping, that is, by choosing $n$ out
of $n$ data points with replacement.

Selecting the points ``without replacement'' instead of ``with
replacement'' is harmless---in fact, it is just a means to avoid
mathematical difficulties induced by the bootstrap; see, for example,
\citet{Ef82,PoRoWo99}.

On the other hand, letting the parameters $a_n$ and $t_n$ depend upon
$n$ offers several degrees of freedom which opens the route for
establishing consistency of the method. To be precise, we will study in
Section~\ref{sec3} the random forest algorithm in two different
regimes. The
first regime is when $t_n<a_n$, which means that trees are not fully
developed. In this case, a proper tuning of $t_n$ ensures the forest's
consistency (Theorem~\ref{theoremconsistencysemidevelopedBRF}).
The second regime occurs when $t_n=a_n$, that is, when trees are fully
grown. In this case, consistency results from an appropriate choice of
the subsample rate $a_n/n$ (Theorem~\ref{Theoremconsistanceforetbreiman}).

So far, we have not made explicit the CART-split criterion used in
Algorithm~\ref{al1}. To properly define it, we let $A$ be a generic cell and
$N_n(A)$ be the number of data points falling in $A$. A cut in $A$ is a
pair $(j,z)$, where $j$ is a dimension in $\{1,\ldots, p\}$ and $z$
is the position of the cut along the $j$th coordinate, within the
limits of $A$. We let $\mathcal{C}_A$ be the set of all such possible
cuts in $A$. Then, with the notation $\mathbf{X}_i = (\mathbf
{X}_i^{(1)},\ldots,
\mathbf{X}_i^{(p)} )$, for any $(j,z) \in\mathcal{C}_A$, the CART-split
criterion [\citet{CART}] takes the form
\begin{eqnarray}
L_{n}(j,z) &= & \frac{1}{N_n(A)} \sum_{i=1}^n
(Y_i - \bar {Y}_{A})^2\mathbh{1}_{\mathbf{X}_i \in A}
\nonumber
\\[-8pt]
\label{definitionempiricalCARTcriterion}
\\[-8pt]
\nonumber
&&{}- \frac{1}{N_n(A)} \sum_{i=1}^n
(Y_i - \bar{Y}_{A_{L}} \mathbh{1}_{\mathbf{X}_i^{(j)} < z} -
\bar{Y}_{A_{R}} \mathbh{1}_{\mathbf
{X}_i^{(j)}
\geq z})^2
\mathbh{1}_{\mathbf{X}_i \in A},
\end{eqnarray}
where $A_L = \{ \mathbf{x}\in A\dvtx  \mathbf{x}^{(j)} < z\}$, $A_R = \{
\mathbf{x}\in A\dvtx  \mathbf{x}^{(j)} \geq z\}$, and $\bar{Y}_{A}$
(resp., $\bar{Y}_{A_{L}}$, $\bar{Y}_{A_{R}}$) is the average of the
$Y_i$'s belonging to $A$ (resp., $A_{L}$, $A_{R}$), with the convention
$0/0=0$. At each cell $A$, the best cut $(j_n^{\star},z_n^{\star})$
is finally selected by maximizing $L_n(j,z)$ over $\mathcal
{M}_{\mathrm{try}}$ and $\mathcal{C}_A$, that is,
\begin{eqnarray*}
&& \bigl(j_n^{\star},z_n^{\star}\bigr) \in
\mathop{\mathop{\operatorname{arg}\operatorname{max}}_{j \in\mathcal{M}_{\mathrm{try}}}}_{(j,z)
\in\mathcal{C}_A}
L_{n}(j,z).
\end{eqnarray*}
To remove ties in the argmax, the best cut is always performed along
the best cut direction $j_n^{\star}$, at the middle of two consecutive
data points.

\section{Main results}\label{sec3}
We consider an additive regression model satisfying the following properties:
\begin{longlist}[(H1)]
\item[(H1)] \textit{The response $Y$ follows
\[
Y = \sum_{j=1}^p m_j\bigl(
\mathbf{X}^{(j)}\bigr) + \varepsilon,
\]
where $\mathbf{X}= (\mathbf{X}^{(1)},\ldots, \mathbf{X}^{(p)})$ is
uniformly
distributed over $[0,1]^p$, $\varepsilon$ is an independent centered
Gaussian noise with finite variance $\sigma^2 >0$ and each component
$m_j$ is continuous.}
\end{longlist}

Additive regression models, which extend linear models, were
popularized by \citet{St85} and \citet{GAM}.
These models, which decompose the regression function as a sum of
univariate functions, are flexible and easy to interpret. They are
acknowledged for providing a good trade-off between model complexity
and calculation time, and accordingly, have been extensively studied
for the last thirty years.
Additive models also play an important role in the context of
high-dimensional data analysis and sparse modeling, where they are
successfully involved in procedures such as the Lasso and various
aggregation schemes; for an overview, see, for example, \citet{HaTiFr09}.
Although random forests fall into the family of nonparametric
procedures, it turns out that the analysis of their properties is
facilitated within the framework of additive models.

Our first result assumes that the total number of leaves $t_n$ in each
tree tends to infinity more slowly than the number of selected data
points $a_n$.

\begin{theorem}\label{theoremconsistencysemidevelopedBRF}
Assume that \textup{(H1)} is satisfied. Then, provided $a_n \to \infty$, $t_n \to \infty$
and $t_n (\log a_n)^9/a_n \to0$, random forests are consistent, that is,
\begin{eqnarray*}
&& \lim_{n \to\infty} \mathbb{E} \bigl[ m_n(\mathbf{X}) - m(
\mathbf {X}) \bigr]^2 = 0.
\end{eqnarray*}
\end{theorem}

It is noteworthy that Theorem~\ref
{theoremconsistencysemidevelopedBRF} still holds with $a_n=n$. In
this case, the subsampling step plays no role in the consistency of the
method. Indeed, controlling the depth of the trees via the parameter
$t_n$ is sufficient to bound the forest error. We note in passing that
an easy adaptation of Theorem~\ref
{theoremconsistencysemidevelopedBRF} shows that the CART algorithm
is consistent under the same assumptions.

The term $(\log a_n)^9$ originates from the Gaussian noise and allows
us to control the noise tail. In the easier situation where the
Gaussian noise is replaced by a bounded random variable, it is easy to
see that the term $(\log a_n)^9$ turns into $\log a_n$, a term which
accounts for the complexity of the tree partition.

Let us now examine the forest behavior in the second regime, where
$t_n=a_n$ (i.e., trees are fully grown), and as before, subsampling is
done at the rate $a_n/n$. The analysis of this regime turns out to be
more complicated, and rests upon assumption (H2) below. We denote
by $Z_i = \mathbh{1}_{\mathbf{X}\stackrel{\Theta}{\leftrightarrow}\mathbf{X}_i}$ the indicator that $\mathbf{X}_i$ falls into the same cell as
$\mathbf{X}$ in
the random tree designed with $\mathcal{D}_n$ and the random parameter
$\Theta$. Similarly, we let $Z_j' = \mathbh{1}_{\mathbf{X}\stackrel
{\Theta
'}{\leftrightarrow} \mathbf{X}_j}$, where $\Theta'$ is an
independent copy
of $\Theta$. Accordingly, we define
\[
\psi_{i,j}(Y_i, Y_j) = \mathbb{E} \bigl[
Z_i Z_j' | \mathbf {X}, \Theta,
\Theta', \mathbf{X}_1,\ldots, \mathbf{X}_n,
Y_i, Y_j \bigr]
\]
and
\[
\psi_{i,j} = \mathbb{E} \bigl[ Z_i
Z_j' | \mathbf{X}, \Theta, \Theta',
\mathbf{X}_1,\ldots, \mathbf{X}_n \bigr].
\]
Finally, for any random variables $W_1$, $W_2$, $Z$, we denote by
$\operatorname{Corr}(W_1$, $ W_2|Z)$ the conditional correlation coefficient
(whenever it exists).

\begin{longlist}[(H2)]
\item[(H2)] \textit{Let $Z_{i,j} = (Z_i, Z_j')$. Then one of the following two conditions
holds\textup{:}}
\begin{longlist}[(H2.1)]
\item[(H2.1)] \textit{One has}
\begin{eqnarray*}
&& \lim_{n \to\infty} (\log a_n)^{2p-2} (\log
n)^{2} \mathbb{E} \Bigl[\mathop{\max_{i,j}}_{i \neq j}
\bigl| \psi_{i,j}(Y_i, Y_j) - \psi_{i,j}\bigr|
\Bigr]^2 = 0.
\end{eqnarray*}
\item[(H2.2)] \textit{There exist a constant $C>0$ and a
sequence $(\gamma_n)_n \to0$ such that, almost surely},
\[
\max_{\ell_1, \ell_2=0,1}\frac{|\operatorname{Corr}(Y_i -
m(\mathbf{X}_i),
\mathbh{1}_{Z_{i,j}=(\ell_1,\ell_2)}| \mathbf{X}_i, \mathbf{X}_j,
Y_j )|}{\mathbb{P}
^{1/2} [Z_{i,j}=(\ell_1,\ell_2) | \mathbf{X}_i, \mathbf{X}_j,
Y_j  ]} \leq\gamma_n
\]
\textit{and}
\begin{eqnarray*}
&& \max_{\ell_1=0,1}\frac{|\operatorname{Corr} ((Y_i - m(\mathbf{X}_i))^2,
\mathbh{1}_{Z_i=\ell_1}| \mathbf{X}_i  )|}{\mathbb{P}^{1/2}
[Z_i =\ell_1
| \mathbf{X}_i  ]} \leq C.
\end{eqnarray*}
\end{longlist}
\end{longlist}

Despite their technical aspect, statements (H2.1) and
(H2.2) have simple interpretations. To understand the meaning of
(H2.1), let us replace the Gaussian noise by a bounded random
variable. A close inspection of Lemma~\ref{lemmaInbis} shows that
(H2.1) may be simply replaced by
\begin{eqnarray*}
&& \lim_{n \to\infty} \mathbb{E} \Bigl[\mathop{\max
_{i,j}}_{i \neq j} \bigl| \psi_{i,j}(Y_i,
Y_j) - \psi_{i,j}\bigr| \Bigr]^2 = 0.
\end{eqnarray*}
Therefore, (H2.1) means that the influence of two $Y$-values on the
probability of connection of two couples of random points tends to zero
as $n\to\infty$.

As for assumption (H2.2), it holds whenever the correlation
between the noise and the probability of connection of two couples of
random points vanishes quickly enough, as $n \to\infty$. Note that,
in the simple case where the partition is independent of the $Y_i$'s,
the correlations in (H2.2) are zero, so that (H2) is
trivially satisfied. This is also verified in the noiseless case, that
is, when $Y=m(\mathbf{X})$.
However, in the most general context, the partitions strongly depend on
the whole sample $\mathcal{D}_n$, and unfortunately, we do not know
whether or not (H2) is satisfied.

\begin{theorem}\label{Theoremconsistanceforetbreiman}
Assume that \textup{(H1)} and \textup{(H2)} are satisfied, and let $t_n=a_n$.
Then, provided $a_n \to\infty$, $t_n\to\infty$ and $a_n \log n /n \to0$, random
forests are consistent, that is,
\[
\lim_{n \to\infty} \mathbb{E} \bigl[ m_n(\mathbf{X}) - m(
\mathbf {X}) \bigr]^2 = 0.
\]
\end{theorem}

To our knowledge, apart from the fact that bootstrapping is replaced by
subsampling, Theorems~\ref{theoremconsistencysemidevelopedBRF} and
\ref{Theoremconsistanceforetbreiman} are the first
consistency results for Breiman's (\citeyear{RF}) forests. Indeed,
most models studied so far are designed independently of $\mathcal
{D}_n$ and
are, consequently, an unrealistic representation of the true procedure.
In fact, understanding Breiman's random forest behavior deserves a more
involved mathematical treatment. Section~\ref{sec4} below offers a thorough
description of the various mathematical forces in action.

Our study also sheds some interesting light on the behavior of forests
when the ambient dimension $p$ is large but the true underlying
dimension of the model is small. To see how, assume that the additive
model (H1) satisfies a sparsity constraint of the form
\begin{eqnarray*}
&& Y = \sum_{j=1}^S m_j
\bigl(\mathbf{X}^{(j)}\bigr) + \varepsilon,
\end{eqnarray*}
where $S<p$ represents the true, but unknown, dimension of the model.
Thus, among the $p$ original features, it is assumed that only the
first (without loss of generality) $S$ variables are informative. Put
differently, $Y$ is assumed to be independent of the last $(p-S)$ variables.
In this dimension reduction context, the ambient dimension $p$ can be
very large, but we believe that the representation is sparse, that is,
that few components of $m$ are nonzero. As such, the value $S$
characterizes the sparsity of the model: the smaller $S$, the sparser $m$.

Proposition~\ref{theoremrelevantvariables} below shows that random
forests nicely adapt to the sparsity setting by asymptotically
performing, with high probability, splits along the $S$ informative variables.

In this proposition, we set $m_{\mathrm{try}}= p$ and, for all
$k$, we denote by $j_{1,n}(\mathbf{X}), \ldots,\break  j_{k,n}(\mathbf
{X})$ the first
$k$ cut directions used to construct the cell containing $\mathbf{X}$, with
the convention that $j_{q,n}(\mathbf{X}) = \infty$ if the cell has
been cut
strictly less than $q$ times.

\begin{proposition} \label{theoremrelevantvariables}
Assume that \textup{(H1)} is satisfied. Let $k \in\mathbb{N}^{\star}$
and $\xi>0$. Assume that there is no interval $[a,b]$ and no $j \in\{
1,\ldots, S\}$ such that $m_j$ is constant on $[a,b]$. Then, with
probability $1-\xi$, for all $n$ large enough, we have, for all $1
\leq q \leq k$,
\begin{eqnarray*}
&& j_{q,n}(\mathbf{X}) \in\{1,\ldots, S\}.
\end{eqnarray*}
\end{proposition}

This proposition provides an interesting perspective on why random
forests are still able to do a good job in a sparse framework.
Since the algorithm selects splits mostly along informative variables,
everything happens as if data were projected onto the vector space
generated by the $S$ informative variables.
Therefore, forests are likely to only depend upon these $S$ variables,
which supports the fact that they have good performance in sparse framework.

It remains that a substantial research effort is still needed to
understand the properties of forests in a high-dimensional setting,
when $p=p_n$ may be substantially larger than the sample size.
Unfortunately, our analysis does not carry over to this context. In
particular, if high-dimensionality is modeled by letting $p_n\to\infty
$, then assumption (H2.1) may be too restrictive since the term
$(\log a_n)^{2p-2}$ will diverge at a fast rate.

%, which, in turn, improves the performance of the method compared to
%other non-adaptive (i.e., whose construction is independent of the
%$Y_i$'s) ones.
%This is in line with the results of \citet{analysis}, who proved that,
%for a simplified model, the performance of the method only depends
%upon the true dimension $S$ and not on the ambient dimension $p$.}

\section{Discussion}\label{sec4}

One of the main difficulties in assessing the mathematical properties
of Breiman's (\citeyear{RF}) forests is that the construction process
of the individual trees strongly depends on both the $X_i$'s and the $Y_i$'s.
For partitions that are independent of the $Y_i$'s, consistency can be
shown by relatively simple means via Stone's (\citeyear{St77})
theorem for local averaging estimates; see also \citet{regression},
Chapter~6. However, our partitions and
trees depend upon the $Y$-values in the data. This makes things
complicated, but mathematically interesting too. Thus, logically, the
proof of Theorem~\ref{Theoremconsistanceforetbreiman} starts with
an adaptation of Stone's (\citeyear{St77}) theorem tailored for
random forests, whereas the proof of Theorem~\ref
{theoremconsistencysemidevelopedBRF} is based on consistency
results of data-dependent partitions developed by \citet{No96}.

Both theorems rely on Proposition~\ref
{variancewithincellforetempirique} below, which stresses an
important feature of the random forest mechanism. It states that the
variation of the regression function $m$ within a cell of a random tree
is small provided $n$ is large enough. To this end, we define, for any
cell $A$, the variation of $m$ within $A$ as
\begin{eqnarray*}
&& \Delta(m,A) = \sup_{\mathbf{x}, \mathbf{x}' \in
A}\bigl|m(\mathbf{x}) - m\bigl(
\mathbf{x}'\bigr)\bigr|.
\end{eqnarray*}
Furthermore, we denote by $A_{n}(\mathbf{X}, \Theta)$ the cell of a tree
built with random parameter $\Theta$ that contains the point $\mathbf{X}$.

\begin{proposition}\label{variancewithincellforetempirique}
Assume that \textup{(H1)} holds. Then, for all $\rho, \xi>0$, there
exists $N \in\mathbb{N}^{\star}$ such that, for all $n >N$,
\begin{eqnarray*}
&& \mathbb{P} \bigl[ \Delta\bigl(m,A_{n}(\mathbf{X}, \Theta)\bigr)
\leq\xi \bigr] \geq1 - \rho.
\end{eqnarray*}
\end{proposition}

It should be noted that in the standard, $Y$-independent analysis of
partitioning regression function estimates, the variance is controlled
by letting the diameters of the tree cells tend to zero in probability.
Instead of such a geometrical assumption, Proposition~\ref
{variancewithincellforetempirique} ensures that the variation of
$m$ inside a cell is small, thereby forcing the approximation error of
the forest to asymptotically approach zero.

While Proposition~\ref{variancewithincellforetempirique} offers a
good control of the approximation error of the forest in both regimes,
a separate analysis is required for the estimation error. In regime $1$
(Theorem~\ref{theoremconsistencysemidevelopedBRF}), the parameter
$t_n$ allows us to control the structure of the tree. This is in line
with standard tree consistency approaches; see, for example, \citet
{DeGyLu96}, Chapter~20. Things are different for the second regime
(Theorem~\ref{Theoremconsistanceforetbreiman}), in which individual
trees are fully grown. In this case, the estimation error is controlled
by forcing the subsampling rate $a_n/n$ to be $\mbox{o}(1/\log n)$,
which is a more unusual requirement and deserves some remarks.

At first, we note that the $\log n$ term in Theorem~\ref
{Theoremconsistanceforetbreiman} is used to control the Gaussian
noise $\varepsilon$.
Thus if the noise is assumed to be a bounded random variable, then the
$\log n$ term disappears, and the condition reduces to $a_n/n \to0$.
The requirement $a_n \log n /n \to0$ guarantees that every single
observation $(\mathbf{X}_i, Y_i)$ is used in the tree construction
with a
probability that becomes small with $n$. It also implies that the query
point $\mathbf{x}$ is not connected to the same data point in a high
proportion of trees. If not, the predicted value at $\mathbf{x}$ would
be influenced too much by one single pair $(\mathbf{X}_i, Y_i)$,
making the
forest inconsistent. In fact, the proof of Theorem~\ref
{Theoremconsistanceforetbreiman} reveals that the estimation error
of a forest estimate is small as soon as the maximum probability of
connection between the query point and all observations is small. Thus
the assumption on the subsampling rate is just a convenient way to
control these probabilities, by ensuring that partitions are dissimilar
enough (i.e., by ensuring that $\mathbf{x}$ is connected with many
data points through the forest).
This idea of diversity among trees was introduced by \citet{RF}, but
is generally difficult to analyze. In our approach, the subsampling is
the key component for imposing tree diversity.

Theorem~\ref{Theoremconsistanceforetbreiman} comes at the price of
assumption (H2), for which we do not know if it is valid in all
generality. On the other hand, Theorem~\ref
{Theoremconsistanceforetbreiman}, which mimics almost perfectly the
algorithm used in practice, is an important step toward understanding
Breiman's random forests. Contrary to most previous works, Theorem~\ref
{Theoremconsistanceforetbreiman} assumes that there is only one
observation per leaf of each individual tree. This implies that the
single trees are eventually not consistent, since standard conditions
for tree consistency require that the number of observations in the
terminal nodes tends to infinity as $n$ grows; see, for example, \citet
{DeGyLu96,regression}. Thus the random forest algorithm aggregates
rough individual tree predictors to build a provably consistent general
architecture.

It is also interesting to note that our results (in particular
Lemma~\ref{consistencyempiricalcutsbestmultiplecutsspurious}) cannot
be directly extended to establish the pointwise consistency of random
forests; that is, for almost all $\mathbf{x}\in[0,1]^d$,
\begin{eqnarray*}
&& \lim_{n \to\infty} \mathbb{E} \bigl[ m_n(\mathbf{x}) -
m(\mathbf{x}) \bigr]^2 = 0.
\end{eqnarray*}
Fixing $\mathbf{x}\in[0,1]^d$, the difficulty results from the fact
that we do not have a control on the diameter of the cell $A_n(\mathbf
{x}, \Theta)$, whereas, since the cells form a partition of $[0,1]^d$,
we have a global control on their diameters. Thus, as highlighted by
\citet{Wa14}, random forests can be inconsistent at some fixed point
$\mathbf{x}\in[0,1]^d$, particularly near the edges, while being
$\mathbb{L}^2$ consistent.

Let us finally mention that all results can be extended to the case
where $\varepsilon$ is a heteroscedastic and sub-Gaussian noise, with
for all $\mathbf{x}\in[0,1]^d$, $\mathbb{V}[\varepsilon| \mathbf
{X}=\mathbf{x}]
\leq\sigma'^2$, for some constant $\sigma'^2$. All proofs can be
readily extended to match this context, at the price of easy technical
adaptations.

%To see this, note that all proofs can be directly extended to that
%case, except that of Lemma
%\ref{Lemmevariancenulleforettheoriquedimensiond} and Lemma
%\ref{CARTconvergencecoupureempiriqueplusieurs}. Simple
%calculations show that Lemma
%\ref{Lemmevariancenulleforettheoriquedimensiond} still holds in
%that case and controlling sub-Gaussian tails allow us to replace
%inequalities (\ref{argument1}),(\ref{argument3}),(\ref{argument3b}),(
%\ref{argument6}) in the proof of Lemma
%\ref{CARTconvergencecoupureempiriqueplusieurs} by similar
%inequalities.

\section{Proof of Theorems~\texorpdfstring{\protect\ref{theoremconsistencysemidevelopedBRF}}{1}
and~\texorpdfstring{\protect\ref{Theoremconsistanceforetbreiman}}{2}}\label{sec5}

For the sake of clarity, proofs of the intermediary results are
gathered in  the supplemental article [\citet{ScBiVE15Supp}]. We
start with some notation.

\subsection{Notation}\label{sec51}
In the sequel, to clarify the notation, we will sometimes write $d =
(d^{(1)},d^{(2)})$ to represent a cut $(j,z)$.

Recall that, for any cell $A$, $\mathcal{C}_{A}$ is the set of all
possible cuts in $A$.
Thus, with this notation, $\mathcal{C}_{[0,1]^p}$ is just the set of
all possible cuts at the root of the tree, that is, all possible
choices $d = (d^{(1)},d^{(2)})$ with $d^{(1)} \in\{1,\ldots, p\}$
and $d^{(2)} \in[0,1]$.

More generally, for any $\mathbf{x}\in[0,1]^p$, we call $\mathcal
{A}_k(\mathbf{x})$ the collection of all possible $k\geq1$
consecutive cuts used to build the cell containing $\mathbf{x}$. Such
a cell is obtained after a sequence of cuts $\mathbf{d}_k = (d_1,\ldots,
d_k)$, where the dependency of $\mathbf{d}_k$ upon $\mathbf{x}$ is
understood. Accordingly, for any $\mathbf{d}_k \in\mathcal{A}_k(\mathbf
{x})$, we let $A(\mathbf{x}, \mathbf{d}_{k})$ be the cell containing
$\mathbf{x}$ built with the particular $k$-tuple of cuts $\mathbf{d}_{k}$.
The proximity between two elements $\mathbf{d}_k$ and $\mathbf{d}'_k$ in
$\mathcal{A}_k(\mathbf{x})$ will be measured via
\begin{eqnarray*}
&& \bigl\|\mathbf{d}_k - \mathbf{d}'_k
\bigr\|_{\infty} = \sup_{1 \leq j \leq k } \max \bigl( \bigl|d_j^{(1)}
- d_j'^{(1)}\bigr|, \bigl|d_j^{(2)}
- d_j'^{(2)}\bigr| \bigr).
\end{eqnarray*}
Accordingly, the distance $d_{\infty}$ between $\mathbf{d}_k\in\mathcal
{A}_k(\mathbf{x})$ and any $\mathcal{A} \subset\mathcal
{A}_k(\mathbf{x})$ is
\begin{eqnarray*}
&& d_{\infty}(\mathbf{d}_k,\mathcal{A}) = \inf
_{\mathbf{z} \in
\mathcal{A}} \|\mathbf{d}_k - \mathbf{z}\|_{\infty}.
\end{eqnarray*}

Remember that $A_{n}(\mathbf{X}, \Theta)$ denotes the cell of a tree
containing $\mathbf{X}$ and designed with random parameter $\Theta$.
Similarly, $A_{k,n}(\mathbf{X}, \Theta)$ is the same cell but where
only the
first $k$ cuts are performed ($k \in\mathbb{N}^{\star}$ is a
parameter to be chosen later). We also denote by $\hat{\mathbf{d}}_{k,n}(\mathbf{X}, \Theta) = (\hat{d}_{1,n}(\mathbf{X}, \Theta
),\ldots,
\hat{d}_{k,n}(\mathbf{X}, \Theta))$ the $k$ cuts used to construct
the cell
$A_{k,n}(\mathbf{X}, \Theta)$.

Recall that, for any cell $A$, the empirical criterion used to split
$A$ in the random forest algorithm is defined in (\ref
{definitionempiricalCARTcriterion}). For any cut $(j,z) \in\mathcal
{C}_A$, we denote the following theoretical version of $L_n(\cdot
,\cdot)$ by
\begin{eqnarray*}
L^{\star}(j,z) & =& \mathbb{V}[ Y | \mathbf{X}\in A ] - \mathbb{P}\bigl[
\mathbf{X}^{(j)} < z | \mathbf{X} \in A \bigr] \mathbb{V}\bigl[Y |
\mathbf{X}^{(j)} < z, \mathbf{X}\in A \bigr]
\\
&&{}- \mathbb{P}\bigl[ \mathbf{X}^{(j)} \geq z | \mathbf{X}\in A \bigr]
\mathbb {V}\bigl[Y | \mathbf{X}^{(j)} \geq z, \mathbf{X}\in A \bigr].
\end{eqnarray*}
Observe that $L^{\star}(\cdot,\cdot)$ does not depend upon the
training set and that, by the strong law of large numbers, $L_n(j,z)
\to L^{\star}(j,z)$ almost surely as $n \to\infty$ for all cuts
$(j,z)\in\mathcal{C}_A$. Therefore, it is natural to define the best
theoretical split $(j^{\star},z^{\star})$ of the cell $A$ as
\begin{eqnarray*}
&& \bigl(j^{\star},z^{\star}\bigr) \in \mathop{\mathop{\operatorname{arg}
\operatorname{min}}_{(j,z) \in \mathcal{C}_A}}_{j \in\mathcal{M}_{\mathrm{try}}} L^{\star}(j,z).
\end{eqnarray*}
In view of this criterion, we define the theoretical random forest as
before, but with consecutive cuts performed by optimizing $L^{\star
}(\cdot,\cdot)$ instead of $L_n(\cdot,\cdot)$. We note that this
new forest does depend on $\Theta$ through $\mathcal{M}_{\mathrm
{try}}$, but not on the
sample $\mathcal D_n$. In particular, the stopping criterion for
dividing cells has to be changed in the theoretical random forest;
instead of stopping when a cell has a single training point, we impose
that each tree of the theoretical forest is stopped at a fixed level $k
\in\mathbb{N}^{\star}$. We also let $A_k^{\star}(\mathbf{X},
\Theta)$ be a
cell of the theoretical random tree at level $k$, containing $\mathbf{X}$,
designed with randomness $\Theta$, and resulting from the $k$
theoretical cuts $\mathbf{d}^{\star}_k(\mathbf{X}, \Theta) =
(d_1^{\star}(\mathbf{X}
, \Theta),\ldots, d_k^{\star}(\mathbf{X}, \Theta))$. Since there
can exist
multiple best cuts at, at least, one node, we call $\mathcal
{A}_k^{\star}(\mathbf{X}, \Theta)$ the set of all $k$-tuples $\mathbf{d}^{\star}_k(\mathbf{X}, \Theta)$ of best theoretical cuts used to build
$A_k^{\star}(\mathbf{X}, \Theta)$.

We are now equipped to prove Proposition~\ref
{variancewithincellforetempirique}.
For reasons of clarity, the proof has been divided
in three steps.
First, we study in Lemma~\ref
{theoremvariancetendvers0forettheorique} the theoretical random forest.
Then we prove in Lemma~\ref
{consistencyempiricalcutsbestmultiplecutsspurious} (via Lemma~\ref
{CARTconvergencecoupureempiriqueplusieurs}) that theoretical
and empirical cuts are close to each other.
Proposition~\ref{variancewithincellforetempirique} is finally
established as a consequence of Lemma~\ref
{theoremvariancetendvers0forettheorique} and
Lemma~\ref{consistencyempiricalcutsbestmultiplecutsspurious}.
Proofs of these lemmas are to be found in the supplemental article
[\citet{ScBiVE15Supp}].

\subsection{Proof of Proposition~\texorpdfstring{\protect\ref{variancewithincellforetempirique}}{2}}

We first need a lemma which states that the variation of $m(\mathbf{X})$
within the cell $A_k^{\star}(\mathbf{X}, \Theta)$ where $\mathbf
{X}$ falls, as
measured by $\Delta(m,A_k^{\star}(\mathbf{X}, \Theta))$, tends to zero.

\begin{lemme} \label{theoremvariancetendvers0forettheorique}
Assume that \textup{(H1)} is satisfied. Then, for all $\mathbf{x}\in[0,1]^p$,
\[
\Delta\bigl(m, A_k^{\star}(\mathbf{x}, \Theta)\bigr) \to0\qquad
\mbox{almost surely, as } k \to\infty.
\]
\end{lemme}

The next step is to show that cuts in theoretical and original forests
are close to each other.
To this end, for any $\mathbf{x}\in[0,1]^p$ and any $k$-tuple of cuts
$\mathbf{d}_k \in\mathcal{A}_k(\mathbf{x})$, we define
\begin{eqnarray*}
L_{n,k}(\mathbf{x}, \mathbf{d}_k) & =& \frac{1}{N_n(A(\mathbf{x}, \mathbf{d}_{k-1}))}
\sum_{i=1}^n (Y_i -
\bar{Y}_{A(\mathbf{x}, \mathbf{d}_{k-1})})^2 \mathbh{1}_{\mathbf{X}_i \in A(\mathbf{x}, \mathbf{d}_{k-1})}
\nonumber
\\
&&{} - \frac{1}{N_n(A(\mathbf{x}, \mathbf{d}_{k-1}))} \sum_{i=1}^n
(Y_i - \bar{Y}_{A_L(\mathbf{x}, \mathbf{d}_{k-1})} \mathbh{1}_{\mathbf{X}
_i^{(d_k^{(1)})} < d_k^{(2)}}
\nonumber
\\
&&\hspace*{40pt}\qquad\qquad\qquad{}- \bar{Y}_{A_R(\mathbf{x}, \mathbf{d}_{k-1})} \mathbh{1}_{\mathbf{X}_i^{(d_k^{(1)})} \geq d_k^{(2)}} )^2
\mathbh{1}_{\mathbf{X}_i
\in A(\mathbf{x}, \mathbf{d}_{k-1})},
\end{eqnarray*}
where $A_L(\mathbf{x}, \mathbf{d}_{k-1}) = A(\mathbf{x}, \mathbf{d}_{k-1})
\cap\{ \mathbf{z}\dvtx  \mathbf{z}^{(d_k^{(1)})} < d_k^{(2)}\}$ and
$A_R(\mathbf{x},
\mathbf{d}_{k-1}) =\break  A(\mathbf{x}, \mathbf{d}_{k-1}) \cap\{ \mathbf{z}\dvtx
\mathbf{z}
^{(d_k^{(1)})} \geq d_k^{(2)}\}$, and where we use the convention
$0/0=0$ when $A(\mathbf{x}, \mathbf{d}_{k-1})$ is empty. Besides, we let
$A(\mathbf{x}, \mathbf{d}_{0}) = [0,1]^p$ in the previous equation. The
quantity $L_{n,k}(\mathbf{x}, \mathbf{d}_k)$ is nothing but the criterion
to maximize in $d_k$ to find the best $k$th cut in the cell $A(\mathbf
{x}, \mathbf{d}_{k-1})$.
Lemma~\ref{CARTconvergencecoupureempiriqueplusieurs} below ensures
that $L_{n,k}(\mathbf{x}, \cdot)$ is stochastically equicontinuous,
for all $\mathbf{x}\in[0,1]^p$.
To this end, for all $\xi>0$, and for all $\mathbf{x}\in[0,1]^p$, we
denote by $\mathcal{A}_{k-1}^{\xi}(\mathbf{x}) \subset\mathcal
{A}_{k-1}(\mathbf{x})$ the set of all $(k-1)$-tuples $\mathbf{d}_{k-1}$
such that the cell $A(\mathbf{x}, \mathbf{d}_{k-1})$ contains a hypercube
of edge length $\xi$. Moreover, we let $\bar{\mathcal{A}}_k^{\xi
}(\mathbf{x}) = \{ \mathbf{d}_k \dvtx  \mathbf{d}_{k-1} \in\mathcal
{A}_{k-1}^{\xi}(\mathbf{x}) \}$ equipped with the norm $\|\mathbf{d}_k\|
_{\infty}$.

\begin{lemme} \label{CARTconvergencecoupureempiriqueplusieurs}
Assume that \textup{(H1)} is satisfied. Fix $\mathbf{x}\in[0,1]^p$, $k
\in\mathbb{N}^{\star}$, and let \mbox{$\xi>0$}.
Then $L_{n,k}(\mathbf{x}, \cdot)$ is stochastically equicontinuous on
$\bar{\mathcal{A}}_k^{\xi}(\mathbf{x})$; that is, for all $\alpha,
\rho> 0$, there exists $\delta>0$ such that
\begin{eqnarray*}
&& \lim_{n \to\infty} \mathbb{P} \Bigl[ \mathop{\sup_{\|\mathbf{d}_k- \mathbf{d}_k'\|_{\infty} \leq\delta}}_{\mathbf{d}_k, \mathbf{d}_k' \in\bar
{\mathcal{A}}_k^{\xi}(\mathbf{x})}
\bigl|L_{n,k}(\mathbf{x}, \mathbf{d}_k) - L_{n,k}\bigl(
\mathbf{x}, \mathbf{d}_k'\bigr) \bigr| > \alpha \Bigr] \leq
\rho.
\end{eqnarray*}
\end{lemme}

Lemma~\ref{CARTconvergencecoupureempiriqueplusieurs} is then used
in Lemma~\ref{consistencyempiricalcutsbestmultiplecutsspurious}
to assess the distance between theoretical and empirical cuts.

\begin{lemme}\label{consistencyempiricalcutsbestmultiplecutsspurious}
Assume that \textup{(H1)} is satisfied. Fix $\xi, \rho>0$ and $k \in
\mathbb{N}^{\star}$. Then there exists $N \in\mathbb{N}^{\star}$
such that, for all $n \geq N$,
\[
\mathbb{P} \bigl[ d_{\infty}\bigl(\hat{\mathbf{d}}_{k,n}(
\mathbf{X}, \Theta), \mathcal{A}_k^{\star}(\mathbf{X}, \Theta)
\bigr) \leq\xi \bigr] \geq1 - \rho.
\]
\end{lemme}

We are now ready to prove Proposition~\ref
{variancewithincellforetempirique}. Fix $\rho, \xi>0$.
Since almost sure convergence implies convergence in probability,
according to Lemma~\ref{theoremvariancetendvers0forettheorique},
there exists $k_0 \in\mathbb{N}^{\star}$ such that
\begin{eqnarray}\label{Finalproof1}
&& \mathbb{P} \bigl[ \Delta\bigl(m, A_{k_0}^{\star}(\mathbf{X},
\Theta)\bigr) \leq\xi \bigr] \geq1 - \rho.
\end{eqnarray}
By Lemma~\ref{consistencyempiricalcutsbestmultiplecutsspurious},
for all $\xi_1>0$, there exists $N \in\mathbb{N}^{\star}$ such
that, for all $n \geq N$,
\begin{eqnarray}\label{Finalproof3}
&& \mathbb{P} \bigl[ d_{\infty}\bigl(\hat{\mathbf{d}}_{k_0, n}(
\mathbf{X}, \Theta), \mathcal {A}_{k_0}^{\star}(\mathbf{X},
\Theta)\bigr) \leq\xi_1 \bigr] \geq1 - \rho.
\end{eqnarray}
Since $m$ is uniformly continuous, we can choose $\xi_1$ sufficiently
small such that, for all $\mathbf{x}\in[0,1]^p$, for all $\mathbf{d}_{k_0}, \mathbf{d}_{k_0}'$
satisfying $d_{\infty}(\mathbf{d}_{k_0}, \mathbf{d}_{k_0}') \leq\xi_1$, we have
\begin{eqnarray}\label{Finalproof4}
&& \bigl\llvert \Delta\bigl(m, A(\mathbf{x}, \mathbf{d}_{k_0})\bigr) -
\Delta\bigl(m, A\bigl(\mathbf{x}, \mathbf{d}_{k_0}'\bigr)
\bigr) \bigr\rrvert \leq\xi.
\end{eqnarray}
Thus, combining inequalities (\ref{Finalproof3}) and (\ref
{Finalproof4}), we obtain
\begin{eqnarray} \label{Finalproof2}
&& \mathbb{P} \bigl[ \bigl\llvert \Delta\bigl(m, A_{k_0,n}(\mathbf{X},
\Theta)\bigr) - \Delta\bigl(m, A_{k_0}^{\star}(\mathbf{X}, \Theta)
\bigr) \bigr\rrvert \leq\xi \bigr] \geq1 - \rho.
\end{eqnarray}
Using the fact that $\Delta(m, A) \leq\Delta(m, A')$ whenever $A
\subset A'$, we deduce from (\ref{Finalproof1}) and (\ref
{Finalproof2}) that, for all $n \geq N$,
\begin{eqnarray*}
&& \mathbb{P} \bigl[ \Delta\bigl(m, A_{n}(\mathbf{X}, \Theta)\bigr)
\leq2 \xi \bigr] \geq1 - 2\rho.
\end{eqnarray*}
This completes the proof of Proposition~\ref
{variancewithincellforetempirique}.

\subsection{Proof of Theorem~\texorpdfstring{\protect\ref{theoremconsistencysemidevelopedBRF}}{1}}

We still need some additional notation. The partition obtained with the
random variable $\Theta$ and the data set $\mathcal{D}_n$ is denoted
by $\mathcal{P}_n(\mathcal{D}_n,\Theta)$, which we abbreviate as
$\mathcal{P}_n(\Theta)$. We let
\begin{eqnarray*}
&& \Pi_n(\Theta) = \bigl\{ \mathcal{P}\bigl( (\mathbf{x}_1,
y_1),\ldots, (\mathbf{x}_n, y_n), \Theta
\bigr)\dvtx (\mathbf{x}_i,y_i) \in[0,1]^d
\times \mathbb{R}\bigr\}
\end{eqnarray*}
be the family of all achievable partitions with random parameter
$\Theta$. Accordingly, we let
\begin{eqnarray*}
&& M\bigl(\Pi_n(\Theta)\bigr) = \max \bigl\lbrace\operatorname{Card}(
\mathcal{P})\dvtx \mathcal{P} \in\Pi_n(\Theta) \bigr\rbrace
\end{eqnarray*}
be the maximal number of terminal nodes among all partitions in $\Pi
_n(\Theta)$. Given a set $\mathbf{z}_1^n = \{\mathbf{z}_1,\ldots,
\mathbf{z}_n\}
\subset[0,1]^d$, $\Gamma(\mathbf{z}_1^n, \Pi_n(\Theta))$ denotes the
number of distinct partitions of $\mathbf{z}_1^n$ induced by elements
of $\Pi
_n(\Theta)$, that is, the number of different partitions $\{ \mathbf{z}_1^n
\cap A\dvtx  A \in\mathcal{P} \}$ of $\mathbf{z}_1^n$, for $\mathcal{P}
\in\Pi
_n(\Theta)$. Consequently, the partitioning number $\Gamma_n(\Pi
_n(\Theta))$ is defined by
\begin{eqnarray*}
&& \Gamma_n\bigl(\Pi_n(\Theta)\bigr) = \max \bigl\lbrace
\Gamma\bigl(\mathbf {z}_1^n, \Pi _n(\Theta)
\bigr)\dvtx \mathbf{z}_1,\ldots, \mathbf{z}_n
\in[0,1]^d \bigr\rbrace.
\end{eqnarray*}
Let $(\beta_n)_n$ be a positive sequence, and define the truncated
operator $T_{\beta_n}$ by
\begin{eqnarray*}
&& \cases{
\displaystyle T_{\beta_n} u = u,  & \quad\mbox{if }
$|u|<\beta_n$,
\vspace*{3pt}\cr
T_{\beta_n} u = \operatorname{sign}(u) \beta_n, &
\quad\mbox{if }$|u|\geq \beta_n$.}
\end{eqnarray*}
Hence\vspace*{1pt} $T_{\beta_n} m_n(\mathbf{X}, \Theta)$, $Y_L = T_L Y$ and
$Y_{i,L} =
T_L Y_i$ are defined unambiguously. We let $\mathcal{F}_n(\Theta)$ be
the set of all functions $f\dvtx  [0,1]^d \to\mathbb{R}$ piecewise
constant on each cell of the partition $\mathcal{P}_n(\Theta)$.
[Notice that $\mathcal{F}_n(\Theta)$ depends on the whole data set.]
Finally, we denote by $\mathcal{I}_{n,\Theta}$ the set of indices of
the data points that are selected during the subsampling step.
Thus the tree estimate $m_n(\mathbf{x},\Theta)$ satisfies
\begin{eqnarray*}
&& m_n( \cdot,\Theta) \in \mathop{\operatorname{arg}\operatorname{min}}_{f \in\mathcal{F}_n(\Theta)}
\frac{1}{a_n} \sum_{i \in\mathcal{I}_{n,\Theta
}} \bigl|f(\mathbf{X}_i)
- Y_i\bigr|^2.
\end{eqnarray*}

The proof of Theorem~\ref{theoremconsistencysemidevelopedBRF} is
based on ideas developed by \citet{No96}, and worked out in Theorem
$10.2$ in \citet{regression}. This theorem, tailored for our context,
is recalled below for the sake of completeness.

\begin{theorem}[{[\citet{regression}]}]
\label{Theorem102gyorfi}
Let $m_n$ and $\mathcal{F}_n(\Theta)$ be as above. Assume that:
\begin{longlist}[(iii)]
\item[(i)] $\lim_{n \to\infty} \beta_n = \infty$;\vspace*{1pt}

\item[(ii)] $\lim_{n \to\infty} \mathbb{E} [\mathop{\inf_{f \in\mathcal{F}_n(\Theta),}}_{\|f\|_{\infty} \leq\beta
_n}\mathbb{E}_{\mathbf{X}}  [f(\mathbf{X}) - m(\mathbf
{X}) ]^2  ] = 0$;\vspace*{1pt}

\item[(iii)] for all $L>0$,
\begin{eqnarray*}
&& \lim_{n \to\infty} \mathbb{E} \biggl[ \mathop{\sup
_{f\in\mathcal{F}_n(\Theta)}}_{\|f\|_{\infty} \leq\beta_n} \biggl| \frac{1}{a_n} \sum
_{i\in\mathcal{I}_{n,\Theta}} \bigl[f( \mathbf{X}_i) - Y_{i,L}
\bigr]^2 - \mathbb{E} \bigl[f(\mathbf{X}) - Y_L
\bigr]^2 \biggr| \biggr] = 0.
\end{eqnarray*}
\end{longlist}
Then
\begin{eqnarray*}
&& \lim_{n \to\infty} \mathbb{E} \bigl[ T_{\beta_n}
m_n(\mathbf{X}, \Theta ) - m(\mathbf{X}) \bigr]^2 = 0.
\end{eqnarray*}
\end{theorem}

Statement (ii) [resp., statement (iii)] allows us to control the
approximation error (resp., the estimation error) of the truncated
estimate. Since the truncated estimate $T_{\beta_n} m_n$ is piecewise
constant on each cell of the partition $\mathcal{P}_n(\Theta)$,
$T_{\beta_n} m_n$ belongs to the set $\mathcal{F}_n(\Theta)$. Thus
the term in (ii) is the classical approximation error.

We are now equipped to prove Theorem~\ref
{theoremconsistencysemidevelopedBRF}. Fix $\xi>0$, and note that
we just have to check statements (i)--(iii) of Theorem~\ref
{Theorem102gyorfi} to prove that the truncated estimate of the
random forest is consistent. Throughout the proof, we let $\beta_n = \|
m\|_{\infty} + \sigma\sqrt{2} (\log a_n)^2$. Clearly, statement
(i) is true.

\subsubsection*{Approximation error}
To prove (ii), let
\[
f_{n,\Theta} = \sum_{A \in\mathcal{P}_n(\Theta)} m(
\mathbf{z}_A) \mathbh{1}_{A},
\]
where $\mathbf{z}_A \in A$ is an arbitrary point picked in cell A. Since,
according to (H1), $\|m\|_{\infty} <\infty$, for all $n$ large
enough such that $\beta_n >\|m\|_{\infty}$, we have
\begin{eqnarray*}
 \mathbb{E}\mathop{\inf_{f \in\mathcal{F}_n(\Theta)}}_{ \|f\|
_{\infty} \leq\beta_n}
\mathbb{E}_{\mathbf{X}} \bigl[f(\mathbf {X}) - m(\mathbf{X})
\bigr]^2  &\leq & \mathbb{E}\mathop{\inf_{f \in\mathcal{F}_n(\Theta)}}_{\|f\|_{\infty} \leq\|m\|_{\infty}}
\mathbb{E}_{\mathbf{X}} \bigl[f(\mathbf{X}) - m(\mathbf{X})
\bigr]^2
\\
& \leq& \mathbb{E} \bigl[f_{\Theta,n}(\mathbf{X}) - m(\mathbf{X})
\bigr]^2
\\
&&{}\bigl(\mbox{since $f_{\Theta,n} \in\mathcal{F}_n(\Theta)$}
\bigr)
\\
&\leq & \mathbb{E} \bigl[m(\mathbf{z}_{A_{n}(\mathbf{X}, \Theta)}) - m(\mathbf{X})
\bigr]^2
\\
& \leq &\mathbb{E} \bigl[ \Delta\bigl(m,A_{n}(\mathbf{X}, \Theta)
\bigr) \bigr]^2
\\
%& \leq\E\big[ \xi^2 \mathbh{1}_{\Delta(m,A_{n}(\bX, \Theta)) \leq
%\xi} + 4\|m\|_{\infty}^2 \mathbh{1}_{\Delta(m,A_{n}(\bX, \Theta)) >
%\xi} \big]\\
& \leq & \xi^2 + 4\|m
\|_{\infty}^2 \mathbb{P} \bigl[ \Delta \bigl(m,A_{n}(
\mathbf{X}, \Theta)\bigr) > \xi \bigr].
\end{eqnarray*}
Thus, using Proposition~\ref{variancewithincellforetempirique}, we
see that for all $n$ large enough,
\begin{eqnarray*}
&& \mathbb{E}\mathop{\inf_{f \in\mathcal{F}_n(\Theta)}}_{\|f\| _{\infty} \leq\beta_n}
\mathbb{E}_{\mathbf{X}} \bigl[f(\mathbf {X}) - m(\mathbf{X})
\bigr]^2 \leq2\xi^2.
\end{eqnarray*}
This establishes (ii).

\subsubsection*{Estimation error}
To prove statement (iii), fix $L>0$. Then, for all $n$ large enough
such that $L <\beta_n$,
\begin{eqnarray*}
&&\mathbb{P}_{\mathbf{X}, \mathcal{D}_n} \biggl( \mathop{\sup_{f \in \mathcal{F}_n(\Theta)}}_{\|f\|_{\infty} \leq\beta_n}
\biggl\llvert \frac{1}{a_n} \sum_{i \in\mathcal{I}_{n,\Theta}} \bigl[f(\mathbf{X}_i) - Y_{i,L} \bigr]^2 - \mathbb{E}
\bigl[f(\mathbf{X}) - Y_L \bigr]^2 \biggr\rrvert > \xi
\biggr)
\\
&&\qquad \leq 8 \exp \biggl[ \log\Gamma_n\bigl(\Pi_n(\Theta)
\bigr) + 2 M\bigl(\Pi _n(\Theta)\bigr) \log \biggl( \frac{333e\beta_n^2}{\xi}
\biggr) - \frac
{a_n \xi^2}{2048\beta_n^4} \biggr]
\\
&& \qquad\quad\mbox{[according to Theorem $9.1$ in \citet{regression}]}
\\
% and Problem $10.4$
&&\qquad \leq8 \exp \biggl[ - \frac{a_n}{\beta_n^4} \biggl( \frac{\xi
^2}{2048} -
\frac{\beta_n^4\log\Gamma_n(\Pi_n)}{a_n} - \frac{2
\beta_n^4 M(\Pi_n)}{a_n} \log \biggl( \frac{333e\beta_n^2}{\xi} \biggr)
\biggr) \biggr].
\end{eqnarray*}
Since each tree has exactly $t_n$ terminal nodes, we have $M(\Pi
_n(\Theta)) = t_n$, and simple calculations show that
\begin{eqnarray*}
&& \Gamma_n\bigl(\Pi_n(\Theta)\bigr) \leq(da_n)^{t_n}.
\end{eqnarray*}
Hence
\begin{eqnarray*}
&&\mathbb{P} \biggl( \mathop{\sup_{f \in\mathcal{F}_n(\Theta)}}_{\|f\|_{\infty} \leq\beta_n}
\biggl\llvert \frac{1}{a_n} \sum_{i \in
\mathcal{I}_{n,\Theta}} \bigl[f(\mathbf{X}_i) - Y_{i,L} \bigr]^2 - \mathbb{E}
\bigl[f(\mathbf{X}) - Y_L \bigr]^2 \biggr\rrvert > \xi
\biggr)
\\
&&\qquad \leq8 \exp \biggl( - \frac{a_nC_{\xi,n}}{\beta_n^4} \biggr),
\end{eqnarray*}
where
\begin{eqnarray*}
C_{\xi,n} & = & \frac{\xi^2}{2048} - 4 \sigma^4
\frac{t_n (\log
(da_n))^9}{a_n} - 8\sigma^4 \frac{t_n(\log a_n)^8}{a_n} \log \biggl(
\frac{666e\sigma^2 (\log a_n)^4}{\xi} \biggr)
\\
& \to & \frac{\xi^2}{2048}\qquad \mbox{as } n \to\infty,
\end{eqnarray*}
by our assumption. Finally, observe that
\begin{eqnarray*}
&& \mathop{\sup_{f \in\mathcal{F}_n(\Theta)}}_{\|f\|
_{\infty} \leq\beta_n} \biggl\llvert
\frac{1}{a_n} \sum_{i \in\mathcal
{I}_{n,\Theta}} \bigl[f(\mathbf{X}_i) - Y_{i,L} \bigr]^2 - \mathbb{E}
\bigl[f(\mathbf{X}) - Y_L \bigr]^2 \biggr\rrvert \leq2(
\beta_n + L)^2,
\end{eqnarray*}
which yields, for all $n$ large enough,
\begin{eqnarray*}
&& \!\!\mathbb{E} \Biggl[ \mathop{\sup_{f \in\mathcal{F}_n(\Theta)}}_{\|f\|_{\infty} \leq\beta_n} \Biggl|
\frac{1}{a_n} \sum_{i=1}^{a_n} \bigl[f(\mathbf{X}_i) - Y_{i,L} \bigr]^2 - \mathbb{E}
\bigl[f(\mathbf{X}) - Y_L \bigr]^2 \Biggr| \Biggr]
\\
&&\!\!\qquad \leq\xi+ 2(\beta_n + L)^2 \mathbb{P} \Biggl[ \mathop{\sup
_{f \in \mathcal{F}_n(\Theta)}}_{\|f\|_{\infty} \leq\beta_n} \Biggl| \frac
{1}{a_n} \sum
_{i=1}^{a_n} \bigl[f(\mathbf{X}_i) -
Y_{i,L} \bigr]^2 - \mathbb{E} \bigl[f(\mathbf{X}) -
Y_L \bigr]^2 \Biggr| > \xi \!\!\Biggr]
\\
&&\!\!\qquad \leq\xi+ 16(\beta_n + L )^2 \exp \biggl( -
\frac{a_nC_{\xi,
n}}{\beta_n^4} \biggr)
\\
&&\!\!\qquad \leq2 \xi.
\end{eqnarray*}
Thus, according to Theorem~\ref{Theorem102gyorfi},
\begin{eqnarray*}
&& \mathbb{E} \bigl[ T_{\beta_n}m_n(\mathbf{X}, \Theta) - m(
\mathbf {X}) \bigr]^2 \to0.
\end{eqnarray*}
\subsubsection*{Untruncated estimate}
It remains to show the consistency of the nontruncated random forest
estimate, and the proof will be complete. For this purpose, note that,
for all $n$ large enough,
\begin{eqnarray*}
\mathbb{E} \bigl[ m_n(\mathbf{X}) - m(\mathbf{X})
\bigr]^2 & =& \mathbb {E} \bigl[ \mathbb{E}_{\Theta
}
\bigl[m_n(\mathbf{X}, \Theta)\bigr] - m(\mathbf{X})
\bigr]^2
\\
& \leq & \mathbb{E} \bigl[ m_n(\mathbf{X}, \Theta) - m(\mathbf{X})
\bigr]^2
\\
&&{} \mbox{(by Jensen's inequality)}
\\
& \leq & \mathbb{E} \bigl[ m_n(\mathbf{X}, \Theta) -
T_{\beta
_n}m_n(\mathbf{X}, \Theta ) \bigr]^2
\\
&&{}+ \mathbb{E} \bigl[ T_{\beta_n}m_n(\mathbf{X}, \Theta) - m(
\mathbf {X}) \bigr]^2
\\
& \leq &\mathbb{E} \bigl[ \bigl[ m_n(\mathbf{X}, \Theta) -
T_{\beta
_n}m_n(\mathbf{X}, \Theta) \bigr]^2
\mathbh{1}_{m_n(\mathbf{X}, \Theta) \geq\beta_n} \bigr] + \xi
\\
& \leq & \mathbb{E} \bigl[ m_n^2(\mathbf{X}, \Theta)
\mathbh {1}_{m_n(\mathbf{X}, \Theta)
\geq\beta_n} \bigr] + \xi
\\
& \leq & \mathbb{E} \bigl[ \mathbb{E} \bigl[ m_n^2(
\mathbf{X}, \Theta) \mathbh{1}_{m_n(\mathbf{X},
\Theta) \geq\beta_n}|\Theta \bigr] \bigr] + \xi.
\end{eqnarray*}
Since $|m_n(\mathbf{X}, \Theta)| \leq\|m\|_{\infty} + \max_{1
\leq i \leq n} |\varepsilon_i|$, we have
\begin{eqnarray*}
&& \mathbb{E} \bigl[ m_n^2(\mathbf{X}, \Theta) \mathbh
{1}_{m_n(\mathbf{X}, \Theta) \geq
\beta_n}|\Theta \bigr]
\\
&&\qquad  \leq\mathbb{E} \Bigl[ \Bigl(2\|m\|_{\infty}^2 + 2\max
_{1 \leq i \leq
a_n} \varepsilon_i^2\Bigr)
\mathbh{1}_{\mathop{\mathop{\max}_{1 \leq i \leq a_n}} \varepsilon_i \geq\sigma\sqrt{2} (\log a_n)^2} \Bigr]
\\
%& \leq2\|m\|_{\infty}^2 \P\big[ \max_{1 \leq i \leq a_n}
%\varepsilon_i \geq\sigma\sqrt{2} (\log a_n)^2 \big] + 2 \E\Big[ (
%\max_{1 \leq i \leq a_n} \varepsilon_i^2) \mathbh{1}_{\max
% _{1 \leq i \leq a_n} \varepsilon_i \geq\sigma\sqrt{2} (\log
%a_n)^2} \Big]\\
%& \leq\E\Big[ \big[ \max_{1 \leq i \leq a_n} \varepsilon_i -
%\sigma\sqrt{2} (\log a_n)^2 \big]^2 \mathbh{1}_{\max_{1 \leq i
%\leq a_n} \varepsilon_i \geq\sigma\sqrt{2} (\log a_n)^2} \Big] + 4
%\varepsilon\\
%& \leq\E\Big[ \big[ \max_{1 \leq i \leq a_n} \varepsilon_i
%\big]^2 \mathbh{1}_{\max_{1 \leq i \leq a_n} \varepsilon_i \geq
%\sigma\sqrt{2} (\log a_n)^2} \Big] + 4\varepsilon\\
&&\qquad
\leq2\|m\|_{\infty}^2 \mathbb{P} \Bigl[ \max
_{1 \leq i \leq a_n} \varepsilon_i \geq\sigma\sqrt{2} (\log
a_n)^2 \Bigr]
\\
&&\qquad\quad {}+ 2 \Bigl(\mathbb{E} \Bigl[ \max_{1 \leq i \leq a_n}
\varepsilon_i^4 \Bigr] \mathbb{P} \Bigl[ \max
_{1 \leq i \leq a_n} \varepsilon_i \geq\sigma\sqrt{2} (\log
a_n)^2 \Bigr] \Bigr)^{1/2}.
\end{eqnarray*}
It is easy to see that
\begin{eqnarray*}
&&  \mathbb{P} \Bigl[ \max_{1 \leq i \leq a_n} \varepsilon_i \geq
\sigma \sqrt{2} (\log a_n)^2 \Bigr] \leq
\frac{a_n^{1 - \log a_n }}{2\sqrt
{\pi} (\log a_n)^2 }.
\end{eqnarray*}
Finally, since the $\varepsilon_i$'s are centered i.i.d. Gaussian
random variables, we have, for all $n$ large enough,
\begin{eqnarray*}
&& \mathbb{E} \bigl[ m_n(\mathbf{X}) - m(\mathbf{X})
\bigr]^2 \\
&&\qquad \leq  \frac{2\|m\|_{\infty}^2
a_n^{1 - \log a_n }}{2\sqrt{\pi} (\log a_n)^2} + \xi + 2 \biggl(3a_n
\sigma^4 \frac{a_n^{1 - \log a_n }}{2\sqrt{\pi}
(\log a_n)^2} \biggr)^{1/2}
\\
&&\qquad  \leq  3\xi.
\end{eqnarray*}
This completes the proof of Theorem~\ref{theoremconsistencysemidevelopedBRF}.

\subsection{Proof of Theorem~\texorpdfstring{\protect\ref{Theoremconsistanceforetbreiman}}{2}}

Recall that each cell contains exactly one data point. Thus, letting
\begin{eqnarray*}
&& W_{ni}(\mathbf{X}) = \mathbb{E}_{\Theta} [ \mathbh
{1}_{\mathbf{X}_i \in
A_n(\mathbf{X}, \Theta)} ],
\end{eqnarray*}
the random forest estimate $m_n$ may be rewritten as
\begin{eqnarray*}
&& m_n(\mathbf{X}) = \sum_{i=1}^n
W_{ni}(\mathbf{X}) Y_i.
\end{eqnarray*}
We have in particular that $\sum_{i=1}^n W_{ni}(\mathbf{X}) = 1$. Thus
\begin{eqnarray*}
\mathbb{E} \bigl[ m_n(\mathbf{X}) - m(\mathbf{X})
\bigr]^2 & \leq &  2 \mathbb{E} \Biggl[ \sum
_{i=1}^n W_{ni}(\mathbf{X})
\bigl(Y_i - m(\mathbf{X}_i) \bigr) \Biggr]^2
\\
&&{}+ 2 \mathbb{E} \Biggl[ \sum_{i=1}^n
W_{ni}(\mathbf{X}) \bigl(m(\mathbf {X}_i) - m(\mathbf{X})
\bigr) \Biggr]^2
\\
& \stackrel{\mathrm{def}}{=} &  2 I_n +2 J_n.
\end{eqnarray*}
\subsubsection*{Approximation error}
Fix $\alpha>0$. To upper bound $J_n$, note that by Jensen's inequality,
\begin{eqnarray*}
J_n &\leq & \mathbb{E} \Biggl[ \sum_{i=1}^n
\mathbh{1}_{\mathbf{X}_i
\in A_n(\mathbf{X},
\Theta)} \bigl(m(\mathbf{X}_i) - m(\mathbf{X})
\bigr)^2 \Biggr]
\\
&\leq & \mathbb{E} \Biggl[ \sum_{i=1}^n
\mathbh{1}_{\mathbf{X}_i \in
A_n(\mathbf{X},
\Theta)} \Delta^2\bigl(m, A_n(
\mathbf{X}, \Theta)\bigr) \Biggr]
\\
&\leq & \mathbb{E} \bigl[ \Delta^2\bigl(m, A_n(\mathbf{X},
\Theta)\bigr) \bigr].
\end{eqnarray*}
So, by definition of $\Delta(m, A_n(\mathbf{X}, \Theta))^2$,
\begin{eqnarray*}
J_n &\leq & 4 \|m\|_{\infty}^2 \mathbb{E}[
\mathbh{1}_{\Delta^2(m,
A_n(\mathbf{X}
, \Theta)) \geq\alpha} ] + \alpha
\\
&\leq & \alpha\bigl(4 \|m\|_{\infty}^2 + 1\bigr),
\end{eqnarray*}
for all $n$ large enough, according to Proposition~\ref
{variancewithincellforetempirique}.
\subsubsection*{Estimation error}
To bound $I_n$ from above, we note that
\begin{eqnarray*}
I_n & = & \mathbb{E} \Biggl[\sum_{i,j=1}^n
W_{ni}(\mathbf{X}) W_{nj}(\mathbf{X}) \bigl(Y_i
- m(\mathbf{X}_i)\bigr) \bigl(Y_j - m(
\mathbf{X}_j)\bigr) \Biggr]
\\
& =&  \mathbb{E} \biggl[\sum_{i=1}
W_{ni}^2(\mathbf{X}) \bigl(Y_i - m(
\mathbf{X}_i)\bigr)^2 \biggr] + I_n',
\end{eqnarray*}
where
\begin{eqnarray*}
&& I_n' = \mathbb{E} \biggl[\mathop{\sum
_{i,j}}_{i \neq j} \mathbh{1}_{\mathbf{X}\stackrel{\Theta}{\leftrightarrow} \mathbf{X}_i}
\mathbh{1}_{\mathbf{X}
\stackrel{\Theta'}{\leftrightarrow} \mathbf{X}_j} \bigl(Y_i - m(\mathbf {X}_i)
\bigr) \bigl(Y_j - m(\mathbf{X}_j)\bigr) \biggr].
\end{eqnarray*}
The term $I'_n$, which involves the double products, is handled
separately in Lemma~\ref{lemmaInbis} below. According to this lemma,
and by assumption (H2), for all $n$ large enough,
\begin{eqnarray*}
&& \bigl|I_n'\bigr| \leq\alpha.
\end{eqnarray*}
Consequently, recalling that $\varepsilon_i = Y_i - m(\mathbf{X}_i)$, we
have, for all $n$ large\vspace*{-2pt} enough,
\begin{eqnarray}
|I_n| & \leq & \alpha+ \mathbb{E} \Biggl[ \sum
_{i=1}^n W_{ni}^2(\mathbf
{X}) \bigl(Y_i - m(\mathbf{X}_i) \bigr)^2
\Biggr]
\nonumber\\[-2pt]
\label{proofprop2eq1}
& \leq & \alpha+ \mathbb{E} \Biggl[ \max_{1\leq\ell\leq n}
W_{n\ell
}(\mathbf{X}) \sum_{i=1}^n
W_{ni}(\mathbf{X}) \varepsilon_i^2 \Biggr]
\\[-2pt]
\nonumber
& \leq & \alpha+ \mathbb{E} \Bigl[ \max_{1\leq\ell\leq n} W_{n\ell
}(
\mathbf{X}) \max_{1\leq i \leq n} \varepsilon_i^2
\Bigr].
\end{eqnarray}
Now, observe that in the subsampling step, there are exactly ${a_n-1\choose  n-1}$ choices to pick a fixed observation $\mathbf{X}_i$. Since
$\mathbf{x}$ and $\mathbf{X}_i$ belong to the same cell only if
$\mathbf{X}_i$ is
selected in the subsampling step, we see\vspace*{-2pt} that
\begin{eqnarray*}
&& \mathbb{P}_{\Theta} [ \mathbf{X}\stackrel{\Theta } {\leftrightarrow}
\mathbf{X}_i ] \leq \frac{{a_n-1 \choose n-1}}{{a_n\choose  n}} = \frac{a_n}{n},
\end{eqnarray*}
where $\mathbb{P}_{\Theta}$ denotes the probability with respect to
$\Theta$, conditional on $\mathbf{X}$ and $\mathcal D_n$.\vspace*{-2pt} So,
\begin{eqnarray}\label{proofprop2eq2}
&& \max_{1 \leq i \leq n} W_{ni}(\mathbf{X})  \leq\max
_{1 \leq i
\leq n } \mathbb{P}_{\Theta} [ \mathbf{X}\stackrel{\Theta }
{\leftrightarrow} \mathbf{X}_i ] \leq\frac{a_n}{n}.
\end{eqnarray}
Thus, combining inequalities (\ref{proofprop2eq1}) and (\ref
{proofprop2eq2}), for all $n$ large enough,
\begin{eqnarray*}
&& |I_n| \leq\alpha+ \frac{a_n}{n} \mathbb{E} \Bigl[ \max
_{1\leq i
\leq n} \varepsilon_i^2 \Bigr].
\end{eqnarray*}
The term inside the brackets is the maximum of $n$ $\chi^2$-squared distributed
random variables. Thus, for some positive constant $C$,
\[
\mathbb{E} \Bigl[ \max_{1\leq i \leq n} \varepsilon_i^2
\Bigr] \leq C\log n;
\]
see, for example, \citet{BoLuMa13}, Chapter~1. We conclude that for all
$n$ large enough,
\[
I_n \leq\alpha+ C \frac{a_n \log n}{n} \leq2 \alpha.
\]
Since $\alpha$ was arbitrary, the proof is complete.

\begin{lemme}\label{lemmaInbis}
Assume that \textup{(H2)} is satisfied. Then, for all $\varepsilon>0$,
and all $n$ large enough, $|I'_n| \leq\alpha$.
\end{lemme}

\begin{pf}%{Proof of Lemma~\protect\ref{lemmaInbis}}
First, assume that (H2.2) is verified.
%To simplify notation, we write $\mathbh{1}_{\bX\stackrel{\Theta}{
%\leftrightarrow} \bX_i}$ instead of $\mathbh{1}_{\bX\in A_n(\bX_i,
%\Theta)}$, keeping in mind the fact that the indicator $\mathbh{1}_{
%\bX\stackrel{\Theta}{\leftrightarrow} \bX_i}$ depends upon the whole
%sample $\mathcal D_n$.
Thus we have for all $\ell_1, \ell_2 \in\{0,1\}$,
\begin{eqnarray*}
&& \operatorname{Corr}\bigl(Y_i - m(\mathbf{X}_i),
\mathbh{1}_{Z_{i,j}=(\ell
_1,\ell
_2)} | \mathbf{X}_i, \mathbf{X}_j,
Y_j\bigr)
\\
&&\qquad  = \frac{\mathbb{E} [(Y_i - m(\mathbf{X}_i)) \mathbh
{1}_{Z_{i,j}=(\ell
_1,\ell_2)}  ]}{\mathbb{V}^{1/2} [Y_i - m(\mathbf{X}_i) |
\mathbf{X}_i, \mathbf{X}_j,
Y_j  ]\mathbb{V}^{1/2} [\mathbh{1}_{Z_{i,j}=(\ell_1,\ell
_2)} | \mathbf{X}
_i, \mathbf{X}_j, Y_j  ]}
\\
&&\qquad = \frac{\mathbb{E} [(Y_i - m(\mathbf{X}_i)) \mathbh
{1}_{Z_{i,j}=(\ell
_1,\ell_2)} | \mathbf{X}_i, \mathbf{X}_j, Y_j  ]}{\sigma
(\mathbb{P}
[Z_{i,j}=(\ell_1,\ell_2) | \mathbf{X}_i, \mathbf{X}_j, Y_j  ]-
\mathbb{P}
[Z_{i,j}=(\ell_1,\ell_2) | \mathbf{X}_i, \mathbf{X}_j, Y_j
]^2 )^{1/2}}
\\
&&\qquad \geq\frac{\mathbb{E} [(Y_i - m(\mathbf{X}_i)) \mathbh
{1}_{Z_{i,j}=(\ell_1,\ell_2)} | \mathbf{X}_i, \mathbf{X}_j, Y_j
 ]}{\sigma\mathbb{P}
^{1/2} [Z_{i,j}=(\ell_1,\ell_2) | \mathbf{X}_i, \mathbf{X}_j,
Y_j  ]},
\end{eqnarray*}
where the first equality comes from the fact that, for all $\ell_1,
\ell_2 \in\{0,1\}$,
\begin{eqnarray*}
&& \operatorname{Cov}\bigl(Y_i - m(\mathbf{X}_i),
\mathbh{1}_{Z_{i,j}=(\ell
_1,\ell
_2)} | \mathbf{X}_i, \mathbf{X}_j,
Y_j\bigr)
\\
&&\qquad  = \mathbb{E} \bigl[\bigl(Y_i - m(\mathbf{X}_i)\bigr)
\mathbh {1}_{Z_{i,j}=(\ell
_1,\ell_2)} | \mathbf{X}_i, \mathbf{X}_j,
Y_j \bigr],
\end{eqnarray*}
since $\mathbb{E}[Y_i - m(\mathbf{X}_i) | \mathbf{X}_i, \mathbf
{X}_j, Y_j] = 0$. Thus, noticing
that, almost surely,
\begin{eqnarray*}
&& \mathbb{E} \bigl[ Y_i - m(\mathbf{X}_i) |
Z_{i,j}, \mathbf {X}_i, \mathbf{X}_j,
Y_j \bigr]
\\
&&\qquad = \sum_{\ell_1, \ell_2 =1}^2 \frac{\mathbb{E} [ (Y_i -
m(\mathbf{X}
_i)) \mathbh{1}_{Z_{i,j}=(\ell_1,\ell_2)} | \mathbf{X}_i, \mathbf
{X}_j, Y_j
 ]}{\mathbb{P} [Z_{i,j} = (\ell_1,\ell_2) | \mathbf{X}_i,
\mathbf{X}_j,
Y_j ]}
\mathbh{1}_{Z_{i,j}=(\ell_1,\ell_2)}
\\
&&\qquad \leq4 \sigma\max_{\ell_1, \ell_2=0,1}\frac{|\operatorname{Corr}
(Y_i - m(\mathbf{X}_i), \mathbh{1}_{Z_{i,j}=(\ell_1,\ell_2)}|
\mathbf{X}_i, \mathbf{X}
_j, Y_j )|}{\mathbb{P}^{1/2} [Z_{i,j}=(\ell_1,\ell_2) | \mathbf
{X}_i, \mathbf{X}_j,
Y_j  ]}
\\
&&\qquad \leq4 \sigma\gamma_n,
\end{eqnarray*}
we conclude that the first statement in (H2.2) implies that,
almost surely,
\begin{eqnarray*}
&&\mathbb{E} \bigl[ Y_i - m(\mathbf{X}_i) |
Z_{i,j}, \mathbf {X}_i, \mathbf{X}_j,
Y_j \bigr] \leq4 \sigma\gamma_n.
\end{eqnarray*}
Similarly, one can prove that the second statement in assumption (H2.2) implies that, almost surely,
\begin{eqnarray*}
&& \mathbb{E} \bigl[ \bigl|Y_i - m(\mathbf{X}_i)\bigr|^2
| \mathbf{X}_i, \mathbh{1}_{\mathbf{X}\stackrel
{\Theta}{\leftrightarrow} \mathbf{X}_i} \bigr] \leq4 C
\sigma^2.
\end{eqnarray*}
Returning to the term $I_n'$, and recalling that $W_{ni}(\mathbf
{X})=\mathbb{E}
_{\Theta}[\mathbh{1}_{\mathbf{X}\stackrel{\Theta}{\leftrightarrow}
\mathbf{X}
_i}]$, we obtain
\begin{eqnarray*}
I_n' & =& \mathbb{E} \biggl[\mathop{\sum
_{i,j}}_{i \neq j} \mathbh {1}_{\mathbf{X}\stackrel{\Theta}{\leftrightarrow} \mathbf{X}_i}
\mathbh{1}_{\mathbf{X}
\stackrel{\Theta'}{\leftrightarrow} \mathbf{X}_j} \bigl(Y_i - m(\mathbf {X}_i)
\bigr) \bigl(Y_j - m(\mathbf{X}_j)\bigr) \biggr]
\\
%& = \sum_{\substack{i,j\\ i \neq j}} \E\left[ W_{ni}(\bX)
%W_{nj}(\bX) (Y_i - m(\bX_i))(Y_j - m(\bX_j)) \right] \\
%& = \sum_{\substack{i,j\\ i \neq j}}\E\left[ \mathbh{1}_{\bX
%\stackrel{\Theta}{\leftrightarrow} \bX_i} \mathbh{1}_{\bX\stackrel{
%\Theta'}{\leftrightarrow} \bX_j} (Y_i - m(\bX_i))(Y_j - m(\bX_j))
%\right] \\
& =& \mathop{
\sum_{i,j}}_{i \neq j}\mathbb{E} \bigl[ \mathbb
{E} \bigl[ \mathbh{1}_{\mathbf{X}\stackrel{\Theta}{\leftrightarrow} \mathbf
{X}_i} \mathbh {1}_{\mathbf{X}\stackrel{\Theta'}{\leftrightarrow} \mathbf{X}_j}
\bigl(Y_i - m(\mathbf{X} _i)\bigr)\\
&&\hspace*{25pt}\quad{}\times \bigl(Y_j
- m(\mathbf{X}_j)\bigr)
 | \mathbf{X}_i, \mathbf{X}_j, Y_i,
\mathbh {1}_{\mathbf{X}\stackrel{\Theta}{\leftrightarrow} \mathbf{X}_i}, \mathbh{1}_{\mathbf{X}
\stackrel{\Theta'}{\leftrightarrow} \mathbf{X}_j} \bigr] \bigr]
\\
& =& \mathop{\sum_{i,j}}_{i \neq j} \mathbb{E}
\bigl[\mathbh {1}_{\mathbf{X}\stackrel{\Theta}{\leftrightarrow} \mathbf{X}_i} \mathbh{1}_{\mathbf{X}
\stackrel{\Theta'}{\leftrightarrow} \mathbf{X}_j} \bigl(Y_i
- m(\mathbf {X}_i)\bigr)
\\
&&\qquad\hspace*{4pt} {}\times\mathbb{E} \bigl[ Y_j - m(\mathbf{X}_j) |
\mathbf{X}_i, \mathbf{X}_j, Y_i,
\mathbh{1}_{\mathbf{X}\stackrel
{\Theta}{\leftrightarrow} \mathbf{X}_i}, \mathbh{1}_{\mathbf
{X}\stackrel{\Theta
'}{\leftrightarrow} \mathbf{X}_j} \bigr] \bigr].
\end{eqnarray*}
Therefore, by assumption (H2.2),
\begin{eqnarray*}
\bigl|I_n'\bigr| & \leq & 4 \sigma\gamma_n \mathop{\sum
_{i,j}}_{i \neq j} \mathbb{E} \bigl[
\mathbh{1}_{\mathbf{X}\stackrel{\Theta
}{\leftrightarrow
} \mathbf{X}_i} \mathbh{1}_{\mathbf{X}\stackrel{\Theta
'}{\leftrightarrow} \mathbf{X}
_j} \bigl|Y_i - m(
\mathbf{X}_i)\bigr| \bigr]
\\
& \leq & \gamma_n \sum_{i=1}^n
\mathbb{E} \bigl[\mathbh{1}_{\mathbf
{X}\stackrel
{\Theta}{\leftrightarrow} \mathbf{X}_i} \bigl|Y_i - m(
\mathbf{X}_i)\bigr| \bigr]
\nonumber
\\
& \leq & \gamma_n \sum_{i=1}^n
\mathbb{E} \bigl[\mathbh{1}_{\mathbf
{X}\stackrel
{\Theta}{\leftrightarrow} \mathbf{X}_i} \mathbb{E} \bigl[ \bigl|Y_i -
m(\mathbf{X}_i)\bigr| | \mathbf{X}_i, \mathbh{1}_{\mathbf{X}\stackrel{\Theta
}{\leftrightarrow} \mathbf{X}_i}
\bigr] \bigr]
\nonumber
\\
& \leq & \gamma_n \sum_{i=1}^n
\mathbb{E} \bigl[\mathbh{1}_{\mathbf
{X}\stackrel
{\Theta}{\leftrightarrow} \mathbf{X}_i} \mathbb{E}^{1/2} \bigl[
\bigl|Y_i - m(\mathbf{X} _i)\bigr|^2 |
\mathbf{X}_i, \mathbh{1}_{\mathbf{X}\stackrel{\Theta
}{\leftrightarrow} \mathbf{X}_i} \bigr] \bigr]
\nonumber
\\
& \leq & 2 \sigma C^{1/2} \gamma_n.
\end{eqnarray*}

This proves the result, provided (H2.2) is true.
Let us now assume that (H2.1) is verified. The key argument is to
note that a data point $\mathbf{X}_i$ can be connected with a random point
$\mathbf{X}$ if $(\mathbf{X}_i, Y_i)$ is selected via the subsampling
procedure and
if there are no other data points in the hyperrectangle defined by
$\mathbf{X}
_i$ and $\mathbf{X}$. Data points $\mathbf{X}_i$ satisfying the
latter geometrical
property are called \textit{layered nearest neighbors} (LNN); see, for
example, \citet{BaSo66}. The connection between LNN and random forests
was first observed by \citet{LiJe06}, and later worked out by
\citet{BiDe10}. It is known, in particular, that the number of LNN
$L_{a_n}(\mathbf{X})$ among $a_n$ data points uniformly distributed on
$[0,1]^d$ satisfies, for some constant $C_1>0$ and for all $n$ large enough,
\begin{eqnarray}
\mathbb{E} \bigl[ L^4_{a_n}(\mathbf{X}) \bigr] & \leq&
a_n \mathbb{P} \bigl[ \mathbf{X}\mathop{\leftrightarrow}\limits^{\Theta}_{\mathrm{LNN}}
\mathbf{X}_j \bigr] + 16 a_n^2 \mathbb {P}
\bigl[ \mathbf{X}\mathop{{\leftrightarrow}}\limits^{\Theta}_{\mathrm{LNN}}
\mathbf{X}_i \bigr] \mathbb{P} \bigl[ \mathbf{X}\mathop{
\leftrightarrow}\limits_{\mathrm{LNN}}^{\Theta} \mathbf{X}_j \bigr]
\nonumber
\\[-8pt]
\label{LNNnumber}
\\[-8pt]
\nonumber
& \leq &  C_1 (\log a_n)^{2d-2};
\end{eqnarray}
see, for example, \citet{BaSo66,BaDeHw05}. Thus we have
\begin{eqnarray*}
&& I_n' %& = \E\left[\sum_{\substack{i,j\\i \neq j}} W_{ni}(\bX) W_{nj}(
%\bX) (Y_i - m(\bX_i))(Y_j - m(\bX_j))\right]\\
%& = \E\left[\sum_{\substack{i,j\\i \neq j}} \mathbh{1}_{\bX
%\stackrel{\Theta}{\leftrightarrow} \bX_i} \mathbh{1}_{\bX\stackrel{
%\Theta'}{\leftrightarrow} \bX_j} (Y_i - m(\bX_i))(Y_j - m(\bX_j))
%\right]\\
= \mathbb{E} \biggl[\mathop{\sum_{i,j}}_{i \neq j}
\mathbh {1}_{\mathbf{X}\stackrel{\Theta}{\leftrightarrow} \mathbf{X}_i} \mathbh{1}_{\mathbf{X}
\stackrel{\Theta'}{\leftrightarrow} \mathbf{X}_j} \mathbh
{1}_{\mathbf{X}_i
\mathop{\leftrightarrow}\limits_{\mathrm{LNN}}^{\Theta}
\mathbf{X}} \mathbh{1}_{\mathbf{X}_j \mathop{\leftrightarrow}\limits_{\mathrm{LNN}}^{\Theta'} \mathbf{X}}\bigl(Y_i - m(
\mathbf{X}_i)\bigr) \bigl(Y_j - m(\mathbf{X}
_j)\bigr) \biggr].
\end{eqnarray*}
Consequently,
\begin{eqnarray*}
I_n' & =& \mathbb{E} \biggl[\mathop{\sum
_{i,j}}_{i \neq j} \bigl(Y_i - m(
\mathbf{X}_i)\bigr) \bigl(Y_j - m(\mathbf{X}_j)
\bigr) \mathbh{1}_{\mathbf{X}_i
\mathop{\leftrightarrow}\limits_{\mathrm{LNN}}^{\Theta} \mathbf{X}} \mathbh {1}_{\mathbf{X}_j
\mathop{\leftrightarrow}\limits_{\mathrm{LNN}}^{\Theta'}
\mathbf{X}}
\\
&&\qquad\hspace*{3pt} {}\times\mathbb{E} \bigl[ \mathbh{1}_{\mathbf{X}\stackrel{\Theta}{\leftrightarrow} \mathbf{X}_i} \mathbh{1}_{\mathbf{X}\stackrel{\Theta
'}{\leftrightarrow} \mathbf{X}_j}
| \mathbf{X}, \Theta, \Theta ', \mathbf{X}_1,\ldots,
\mathbf{X}_n, Y_i, Y_j \bigr] \biggr],
\end{eqnarray*}
where\vspace*{1.5pt} $\mathbf{X}_i \mathop{\leftrightarrow}\limits_{\mathrm{LNN}}^{\Theta} \mathbf{X}$ is the event where $\mathbf{X}_i$ is
selected by
the subsampling and is also a LNN of $\mathbf{X}$. Next, with the
notation of
assumption\vspace*{-2pt} (H2),
\begin{eqnarray*}
I_n' & =& \mathbb{E} \biggl[\mathop{\sum
_{i,j}}_{i \neq j} \bigl(Y_i - m(
\mathbf{X}_i)\bigr) \bigl(Y_j - m(\mathbf{X}_j)
\bigr) \mathbh{1}_{\mathbf{X}_i
\mathop{\leftrightarrow}\limits_{\mathrm{LNN}}^{\Theta} \mathbf{X}} \mathbh {1}_{\mathbf{X}_j
\mathop{\leftrightarrow}\limits_{\mathrm{LNN}}^{\Theta'}
\mathbf{X}}
\psi_{i,j}(Y_i, Y_j) \biggr]
\\[-2pt]
& = & \mathbb{E} \biggl[\mathop{\sum_{i,j}}_{i \neq j}
\bigl(Y_i - m(\mathbf{X} _i)\bigr)
\bigl(Y_j - m(\mathbf{X}_j)\bigr) \mathbh{1}_{\mathbf{X}_i
\mathop{\leftrightarrow}\limits_{\mathrm{LNN}}^{\Theta} \mathbf{X}}
\mathbh {1}_{\mathbf{X}_j
\mathop{\leftrightarrow}\limits_{\mathrm{LNN}}^{\Theta'}
\mathbf{X}} \psi_{i,j} \biggr]
\\[-2pt]
&&{}+ \mathbb{E} \biggl[\mathop{\sum_{i,j}}_{i \neq j}
\bigl(Y_i - m(\mathbf{X}_i)\bigr) \bigl(Y_j
- m(\mathbf{X}_j)\bigr) \mathbh{1}_{\mathbf{X}_i
\mathop{\leftrightarrow}\limits_{\mathrm{LNN}}^{\Theta} \mathbf{X}} \mathbh
{1}_{\mathbf{X}_j
\mathop{\leftrightarrow}\limits_{\mathrm{LNN}}^{\Theta'}
\mathbf{X}}
\bigl( \psi_{i,j}(Y_i, Y_j) -
\psi_{i,j}\bigr) \biggr].
\end{eqnarray*}
The first term is easily seen to be zero\vspace*{-2pt} since
\begin{eqnarray*}
&& \mathbb{E} \biggl[\mathop{\sum_{i,j}}_{i \neq j}
\bigl(Y_i - m(\mathbf{X} _i)\bigr)
\bigl(Y_j - m(\mathbf{X}_j)\bigr) \mathbh{1}_{\mathbf{X}_i
\mathop{\leftrightarrow}\limits_{\mathrm{LNN}}^{\Theta} \mathbf{X}}
\mathbh {1}_{\mathbf{X}_j
\mathop{\leftrightarrow}\limits_{\mathrm{LNN}}^{\Theta'}
\mathbf{X}} \psi\bigl(\mathbf{X}, \Theta, \Theta',
\mathbf{X}_1,\ldots, \mathbf{X}_n\bigr) \biggr]
\\[-2pt]
&&\qquad = \mathop{\sum_{i,j}}_{i \neq j} \mathbb{E}
\bigl[ \mathbh {1}_{\mathbf{X}_i
\mathop{\leftrightarrow}\limits_{\mathrm{LNN}}^{\Theta}
\mathbf{X}} \mathbh{1}_{\mathbf{X}_j
\mathop{\leftrightarrow}\limits_{\mathrm{LNN}}^{\Theta'} \mathbf{X}} \psi_{i,j}
\\[-2pt]
&&\qquad\qquad\hspace*{13pt}{}\times\mathbb{E} \bigl[ \bigl(Y_i - m(\mathbf{X}_i)
\bigr) \bigl(Y_j - m(\mathbf{X}_j)\bigr) | \mathbf{X},
\mathbf{X}_1,\ldots, \mathbf{X}_n, \Theta,
\Theta' \bigr] \bigr]
\\
&&\qquad = 0.
\end{eqnarray*}
Therefore,
\begin{eqnarray*}
\bigl|I_n'\bigr| & \leq & \mathbb{E} \biggl[\mathop{\sum
_{i,j}}_{i \neq j} \bigl|Y_i - m(
\mathbf{X}_i)\bigr|\bigl|Y_j - m(\mathbf{X}_j)\bigr|
\mathbh{1}_{\mathbf{X}_i
\mathop{\leftrightarrow}\limits_{\mathrm{LNN}}^{\Theta} \mathbf{X}} \mathbh {1}_{\mathbf{X}_j \mathop{\leftrightarrow}\limits_{\mathrm{LNN}}^{\Theta'} \mathbf{X}}
\bigl| \psi_{i,j}(Y_i, Y_j) -
\psi_{i,j}\bigr| \biggr]
\\
& \leq & \mathbb{E} \biggl[ \max_{1 \leq\ell\leq n}\bigl|Y_i - m(
\mathbf{X}_i)\bigr|^2 \mathop{\max_{i,j}}_{i \neq j}
\bigl| \psi_{i,j}( Y_i, Y_j) - \psi_{i,j}\bigr|
%\\
%&&{}\times
\mathop{\sum_{i,j}}_{i \neq j}
\mathbh{1}_{\mathbf{X}_i \mathop{\leftrightarrow}\limits_{\mathrm{LNN}}^{\Theta} \mathbf{X}} \mathbh{1}_{\mathbf{X}_j
\mathop{\leftrightarrow}\limits_{\mathrm{LNN}}^{\Theta'} \mathbf{X}} \biggr].
\end{eqnarray*}
Now, observe that
\begin{eqnarray*}
&& \mathop{\sum_{i,j}}_{i \neq j}
\mathbh{1}_{\mathbf{X}_i \mathop{\leftrightarrow}\limits_{\mathrm{LNN}}^{\Theta} \mathbf{X}} \mathbh{1}_{\mathbf{X}_j \mathop{\leftrightarrow}\limits_{\mathrm{LNN}}^{\Theta'} \mathbf{X}} \leq L^2_{a_n}(
\mathbf{X}).
\end{eqnarray*}
Consequently,
\begin{eqnarray}
\bigl|I_n'\bigr| & \leq & \mathbb{E}^{1/2} \Bigl[
L^4_{a_n}(\mathbf{X}) \max_{1
\leq
\ell\leq n}\bigl|Y_i
- m(\mathbf{X}_i)\bigr|^4 \Bigr]
\nonumber
\\[-8pt]
\label{equationtemporaire}
\\[-8pt]
\nonumber
&{}& \times\mathbb{E}^{1/2} \Bigl[\mathop{\max_{i,j}}_{i \neq j}
\bigl| \psi_{i,j}( Y_i, Y_j) - \psi_{i,j}\bigr|
\Bigr]^2.
\end{eqnarray}
Simple calculations reveal that there exists $C_1>0$ such that, for all $n$,
\begin{eqnarray} \label{LNNnumber2}
&& \mathbb{E} \Bigl[ \max_{1 \leq\ell\leq n}\bigl|Y_i - m(
\mathbf{X}_i)\bigr|^4 \Bigr] \leq C_1 (\log
n)^2.
\end{eqnarray}
Thus, by inequalities (\ref{LNNnumber}) and (\ref{LNNnumber2}),
the first term in (\ref{equationtemporaire}) can be upper bounded as follows:
\begin{eqnarray*}
&& \mathbb{E}^{1/2} \Bigl[ L^4_{a_n}(\mathbf{X})
\max_{1 \leq\ell
\leq
n}\bigl|Y_i - m(\mathbf{X}_i)\bigr|^4
\Bigr]
\\
&&\qquad = \mathbb{E}^{1/2} \Bigl[ L^4_{a_n}(\mathbf{X})
\mathbb{E} \Bigl[ \max_{1 \leq
\ell\leq n} \bigl|Y_i - m(
\mathbf{X}_i)\bigr|^4 |\mathbf{X}, \mathbf {X}_1,\ldots, \mathbf{X}_n \Bigr] \Bigr]
\\
%& \leq(\log n)^{\alpha/2} \E^{1/2}\Bigg[ L_{a_n}(\bX)^4 \Bigg]\\
%& \leq(\log n)^{\alpha/2} \bigg( a_n \P\big[ \bX\underset{{\tiny
%\mbox{LNN}}}{\stackrel{\Theta}{\leftrightarrow}} \bX_j \big] + 16 a_n^2
%\P\big[ \bX\underset{{\tiny\mbox{LNN}}}{\stackrel{\Theta}{
%\leftrightarrow}} \bX_i \big] \P\big[ \bX\underset{{\tiny
%\mbox{LNN}}}{\stackrel{\Theta}{\leftrightarrow}} \bX_j \big]
%\bigg)^{1/2} \\
&&\qquad
\leq C' (\log n) (\log a_n)^{d-1}.
\end{eqnarray*}
Finally,
\begin{eqnarray*}
&&\bigl|I_n'\bigr| \leq C' (\log
a_n)^{d-1} (\log n)^{\alpha/2} \mathbb
{E}^{1/2} \Bigl[\mathop{\max_{i,j}}_{i \neq j} \bigl|
\psi_{i,j}(Y_i, Y_j) - \psi_{i,j}\bigr|
\Bigr]^2,
\end{eqnarray*}
which tends to zero by assumption.
\end{pf}

\section*{Acknowledgments}
We greatly thank two referees for valuable comments
and insightful suggestions.

%\begin{supplement}[id=suppA]
%\stitle{Technical results}
%\slink[doi]{COMPLETED BY THE TYPESETTER}
%\sdescription{Proofs of technical results}
%\end{supplement}

\begin{supplement}[id=suppA]
\stitle{Supplement to ``Consistency of random forests''\\}
\slink[doi]{10.1214/15-AOS1321SUPP} %[doi,text={...}] - jei reikia
%suskaldyti doi
\sdatatype{.pdf}
\sfilename{aos1321\_supp.pdf}
\sdescription{Proofs of technical results.}
\end{supplement}

% imsref loaded by daiva.urboniene, 2015-03-04 13:27:24
% imsref loaded by daiva.urboniene, 2015-03-09 11:13:51

%\begin{appendix}
%\section{}
%\end{appendix}

%\section*{Acknowledgments}

%\begin{thebibliography}{99}
%\bibitem{r1}
%\bibitem{r1}
%\end{thebibliography}

\printaddresses
\end{document}